\documentclass[12pt]{article}

        \usepackage{amsfonts}
        \usepackage{amssymb}

\newtheorem{thm}{Theorem}
\newtheorem{cor}[thm]{Corollary}
\newtheorem{lemma}[thm]{Lemma}

\newtheorem{prop}[thm]{Proposition}

        \def\pf{\medbreak\noindent{\bf Proof:}\enspace}
        
        \def\qed{{\bf QED}}
        \def\iff{\Longleftrightarrow}
        \def\half{{\textstyle \frac{1}{2}}}
        
        \def\fm{\lfloor m \rfloor}
   
\def\dffx{\frac{d~}{dx}}
\def\dffy{\frac{d~}{dy}}
     
\def\av{{\rm av}}
  \def\raw{\rightarrow}

 \def\ds{\displaystyle}
 \def\fp{{\textstyle \frac{1}{p}}}
 \def\ts{\textstyle}

        \title{Study of a Class of Regularizations \\ of
   $1/|x|$ using Gaussian Integrals}

        \author {Mary Beth Ruskai \thanks{supported by National Science
        Foundation Grant DMS-97-06981} \\ Department of
        Mathematics \\ University of Massachusetts  Lowell \\ Lowell,
        MA  01854 USA \\ {\normalsize bruskai@cs.uml.edu}
        \and Elisabeth Werner \\Department of Mathematics \\
    Case Western Reserve University \\
    Cleveland, OH 44106 USA \\
    {\normalsize emw2@po.cwru.edu}\\
    and\\
    Universit\'{e} de Lille 1\\
    Ufr de Math\'{e}matique\\
    59655 Villeneuve d'Ascq, France}

\date{16 February 1999 \\ revised 9 September 1999 and 15 November 1999}

\begin{document}

 \maketitle

\begin{abstract}
 This paper presents a comprehensive study of the functions
$ V_m^p(x) = \frac{pe^{x^p}}{\Gamma(m+1)}
            \int_x^\infty  (t^p-x^p)^me^{-t^p} dt $ for $x > 0$,
$m > -1$ and $p > 0$.  For large $x$ these functions approximate
$x^{1-p}$.  The case $p=2$ is of particular importance because
the functions  $V_m^2(x) \approx  1/x$ can be regarded as
one-dimensional regularizations of the Coulomb potential $1/|x|$
which are finite at the origin for $m > - \half$.

   The limiting behavior and monotonicity properties of these functions
are discussed in terms of their dependence on $m$ and $p$ as well as $x$.
Several classes of inequalities, some of which provide tight bounds,
are established.   Some differential equations and recursion relations
satisfied by these functions are given.  The recursion relations give
rise to two classes of polynomials, one of which is related to confluent
hypergeometric functions.   Finally, it is shown that, for integer $m$,
the function
$1/V_m^2(x)$ is convex in $x$ and this implies an analogue of the
triangle inequality.   Some comments are made about the range of $p$
and $m$ to which this convexity result can be extended and several
related questions are raised.

\end{abstract}

 \tableofcontents

 \pagebreak

\section{Introduction}

\subsection{Definitions and background}

In this paper we study the functions
\begin{eqnarray}\label{vmequiv1}
 V_m(x) & = & \frac{2e^{x^2}}{\Gamma(m+1)}
            \int_x^\infty  (t^2-x^2)^me^{-t^2} dt,
     ~~~~ m > -1 \\
 V_{-1}(x) & = &    \frac{1}{|x|}   \nonumber
\end{eqnarray}
and their generalizations,
\begin{eqnarray}\label{vmp1}
   V_m^p(x) & = & \frac{pe^{x^p}}{\Gamma(m+1)}
            \int_x^\infty  (t^p-x^p)^me^{-t^p} dt  \\
       V_{-1}^p(x) & = &   x^{1-p},  \nonumber
\end{eqnarray}
for $0 < p < \infty$.
These functions are well-defined for $x > 0$ and can be
extended to complex $m$ with $\Re (m) > - 1 $.  For $\Re (m) > - \half $
they are also well-defined for $x = 0$.  Using symmetry or the
equivalent forms  (\ref{vmequiv2}) and (\ref{vmp2}) below, they
can be extended to even functions on ${\bf R}$
or ${\bf R} \setminus \{0\}$.     However, it suffices to consider only
non-negative $x$ and in this paper we restrict ourselves to that.
We also restrict ourselves to real $m$.

Letting $p=2$ in  (\ref{vmp1}) yields (\ref{vmequiv1}).
However, because this case is more important in applications
we often drop the superscript and and simply write
 $V_m(x)$ for $V_m^2(x)$.

Our interest was motivated by studies  of atoms in magnetic fields
where these functions arise naturally for integer $m$.
$V_m$ can be regarded as a (two-dimensional) expectation of the
(three-dimensional) Coulomb potential $1/|{\bf r}|$ with the state
$\gamma_m(r,\theta) = \frac{1}{\sqrt{ \pi m!}} e^{-im \theta}r^{m}
e^{-r^2/2}$  (where we have used cylindrical coordinates ${\bf r} =
(x,r,\theta)$ with the non-standard convention $r = \sqrt{y^2 + z^2}$
if ${\bf r} = (x,y,z)$ in rectangular coordinates).
The state $\gamma_m$ describes an electron  in the lowest of the
so-called ``Landau levels'' with angular momentum
$m$ in the direction of the field.
In this context, it is natural to rewrite (\ref{vmequiv1}) in the form
\begin{eqnarray}
   V_m(x) & = &\frac{2}{\Gamma(m+1)}\int_0^\infty
        \frac{r^{2m} e^{-r^2} }{\sqrt{x^2 + r^2}} r dr,  \label{vmdef}  \\
     & = &  \frac{1}{\Gamma(m+1)} \int_0^\infty
          \frac{u^m e^{-u} }{\sqrt{x^2 + u}} du  \label{vmequiv2}
\end{eqnarray}
for $m > -1$.  In this form,  it is
easy to see that  $ V_m(x) \approx 1/|x|$ for large $x$.
The importance of $V_m$  goes back
at least to Schiff and Snyder \cite{SS} and played an essential
role in the Avron, Herbst and Simon \cite{AHS} study of the
energy asymptotics of hydrogen in a strong magnetic field.
More recently work in astrophysics and the work of Lieb,
Solovej and Yngvason \cite{LSY} on asymptotics of many-electron atoms
in strong magnetic fields has renewed interest in this subject.
Motivated by the LSY work, Brummelhuis and Ruskai \cite{BR1,BR2}
have developed one-dimensional models of many-electron atoms in
strong magnetic fields using the functions $V_m(x)$ as one-dimensional
analogues of the Coulomb potential.

In the case of many-electron atoms, the anti-symmetry required by
the Pauli exclusion principle suggests replacing the simple
``one-electron'' expectation above by an N-electron analogue in
which the state $\gamma_m$ is replaced by a Slater determinant
of such states.  This is discussed in detail in \cite{BR2} where
it is shown that, in the simple case corresponding to
$m = 0 \ldots {N-1}$, the analogous one-dimensional potentials
have the form
\begin{eqnarray}\label{def.veff}
V_{\av}^N(x) = \frac{1}{N} \sum_{m=0}^{N-1} V_m(x).
\end{eqnarray}
In Section \ref{sect:recurs} we obtain recursion relations for
$V_m$  which, in addition  to being of considerable interest in
their own right, are extremely useful for studying potentials
of the form (\ref{def.veff}).

For $m=0$ the function
\begin{eqnarray}\label{mills}
  \frac{1}{\sqrt{2}} V_0\left( \frac{x}{\sqrt{2}}\right) = e^{x^2/2}
\int_x^\infty e^{-t^2/2} dt
\end{eqnarray}
occurs in many other contexts and is sometimes called the
``Mills ratio'' \cite{M}.   Although it has been extensively
studied, the class of inequalities we consider in Section \ref{sect:V0}
appears to be new (although some of our bounds coincide with
known inequalities in other classes) and the realization
that $1/V_0(x)$ is convex seems to be relatively recent
\cite{SzW,W,BR1}.

The replacement of $x^2$ by $x^p$ in (\ref{mills}) has been considered
by Gautschi \cite{G} and  Mascioni \cite{Mas}, who (after seeing
the preprint \cite {RuW}) extended the
results of Section \ref{sect:V0} to this situation.

For the analysis of this generalization, it is useful to observe that
(\ref{vmp1}) can  be rewritten as
\begin{eqnarray}
   V_m^p(x) = \frac{1}{\Gamma(m+1)} \int_0^\infty
          \frac{u^m e^{-u} }{(x^p + u)^{\frac{p-1}{p}}} du  \label{vmp2}
\end{eqnarray}
In this form, it is easy to verify that $V_m^p(x) \approx x^{p-1}$
for large $x$ and that for $p = 1$, $V_m^1(x) = 1$ for all $m$ .

\bigskip

Our first result shows that $V_m^p(x)$ is continuous in $m$
and that our definition for $m = -1$ is natural.
\begin{prop}
For all $x > 0$,
$\displaystyle{\lim_{m \rightarrow -1^+} x^{p-1} V_m^p(x) = 1}$.
\end{prop}
\pf Note that $\ds{~
   \frac{x^{p-1}} { (x^p + u)^{\frac{p-1}{p}} } =
    \frac{1}{\left( 1 + \frac{u}{x^p} \right)^{\frac{p-1}{p}}} ~~}$
so that (\ref{vmp2}) implies
\begin{eqnarray}\label{m1p.lim.pf.eq1}
  1 - x^{p-1}  V_m^p(x) = \frac{1}{\Gamma(m+1)} \int_0^\infty u^m e^{-u}
\left[ 1- \left( 1 + \frac{u}{x^p} \right)^{\frac{1-p}{p}} \right] du.
\end{eqnarray}
Since $\Gamma(m+1)$ becomes infinite as $m \raw -1$, the
desired result follows if the integral on the right above remains finite.
 To see that this is true, it is convenient
to let $g(z) = \frac{1}{z}\left(1 - (1+z)^{\frac{1-p}{p}}\right)$
and note that (\ref{m1p.lim.pf.eq1}) implies
\begin{eqnarray}\label{m1p.lim.pf.eq2}
\left| 1 - x^{p-1}  V_m^p(x) \right|
 \leq   \frac{1}{\Gamma(m+1)} \int_0^\infty \frac{u^{m+1} e^{-u}}{x^p}
  ~ \left| g\left(  \frac{u}{x^p}  \right) \right| du
\end{eqnarray}
It is easy to see that the large $u$ portion of this integral causes
no problems since $|g(z)|$ is bounded by a polynomial in $z$ when
$z > 1$. (For $p \geq 1$, it is bounded uniformly by $1$;
for $0 < p < 1$,
it is bounded by a polynomial, namely
$|g(z)| \leq 1+(1+z)^k$, where $k \in {\bf
N}$ , $k \geq \frac{1}{p}$.)
To see that it is also well-behaved for small $u$, we
first note that for $p > 0$,
$(1+z)^{\frac{1-p}{p}}$ is analytic for $\Re(z)  > -1$.  Then
$g(z)$ has a removable
singularity at $z = 0$ and can be extended to an analytic function on
$\Re(z) > -1$.  Thus, for small $u$, the integrand behaves
like $x^{-p} u^{m+1} e^{-u}$ which ensures that the integral in
(\ref{m1p.lim.pf.eq2}) is finite for $m = -1$.
 \qed

\bigskip
The rest of this paper is organized as follows.  In the next part
of this section we summarize the properties of $V_m$ in the
important case $p = 2$.  We then conclude the Introduction with a
summary of convexity results, including some open questions. In
Section \ref{sect:gen.p} we state and prove the basic properties
of $ V_m^p$ for general $p$.
 In Section \ref{sect:recurs} we derive recursion relations for
$V_m^p$ and study their consequences.   Among these is a
connection with confluent hypergeometric functions.  In Section
\ref{sect:V0} we prove some optimal bounds for $V_0$. The optimal upper
bound had been established earlier independently by Wirth \cite{W}
and by  Szarek and Werner \cite{SzW}  who also showed that the upper bound
is equivalent to  the convexity of $1/V_0$.
In Section  \ref{sect:rat.bnd} we discuss several classes of inequalities,
beginning with optimal bounds on $V_0(x)$.
We then consider optimal bounds on the ratio
$R_{m}(x)=V_{m}(x) / V_{m-1}(x)$
and show that these have important consequences.  In particular,
we show that the upper bound is equivalent to the convexity (in
$x$) of $1/V_m(x)$ and that the ratios increase with $x$.  Proofs of the
ratio bounds are then given in Section \ref{sect:rat.bnd.pfs} where we
also consider extensions to other $p$.  Because the proof of the ratio
bounds is via induction on $m$, the results of Sections 6 and 7 are only
established for integer $m$.  However, we believe that they hold for all
$m > -1$.

\subsection{Properties of $V_m(x)$}\label{subsect:propvm}

We now summarize some properties of $V_m(x)$ along with comments
about the history and brief remarks about the proofs.
Unless otherwise stated, these properties hold for $m > -1$ and $x > 0$.

\begin{itemize}

  \item[a)] $~~\ds{\frac{1}{\sqrt{x^2 + m}}  > V_m(x)   > \frac{1}
{\sqrt{x^2 + m + 1}}}$
\vspace{0.2cm}
\newline where the first inequality holds for  $m > 0$ and the
second for $m > -1.$

To prove the upper bound, which appears to be new, observe that
$\mu=[u^{m-1} e^{-u}/\Gamma(m)]du$ is a probability measure on $(0, \infty)$.
For fixed $x$, one can then apply Jensen's inequality to
the concave function $f_x(u) = u (u + x^2)^{-1/2}$ to obtain
$$V_m(x)=\frac{1}{m}\int_0^\infty f_x(u) d\mu(u)
\leq \frac{1}{m}f_x \biggl(\int_0^\infty
\frac{u^me^{-u}du}{\Gamma(m)}\biggr)$$
$$=\frac{1}{m}f_x(m)=\frac{1}{\sqrt{x^2+m}}.$$

The lower bound was proved earlier (at least for integer $m$)
by Avron, Herbst and Simon \cite{AHS}
who applied a similar argument to the probability measure
$[u^{m} e^{-u}/\Gamma(m+1)]du$ and the convex function
$f_x(u) = (u +x^2)^{-1/2}$.

\item[b)] $V_m(x)$ is decreasing in $m$.  In particular,
 $\ds{ V_{m+1}(x) < V_m(x) < \frac{1}{x} }$.

The first inequality follows easily from property (a) which implies
\linebreak $V_m(x) < \frac{1}{\sqrt{x^2+m}} < V_{m-1}(x).$
Alternatively, one could use integration
by parts on (\ref{vmequiv2}).
The second inequality is easily
verified from the integral representation (\ref{vmequiv2}).
That $V_m(x)$ also decreases with $m$ for non-integer jumps is more
difficult, and the proof is postponed to Section \ref{sect:gen.p}
where it follows from the more general Theorem \ref{thm:vmp.mono}.

   \item[c)] The expression $ m V_m(x)$ is increasing in $m > -1$, $m
\in{\bf R}$.
\newline
For integer jumps this holds for $m \geq -1$. Indeed,
it is obvious that $-V_{-1} < 0 \cdot V{_0} < V_1$.
For integer jumps with $m \geq 1$, one can use property (a)
to see that
\begin{eqnarray}
  m V_m(x) > \frac{m}{\sqrt{x^2 + m+ 1}} >
    \frac{m-1}{\sqrt{x^2 + m - 1}} > (m-1) V_{m-1}(x).
\end{eqnarray}
The  proof for general $m$ is postponed to Theorem \ref{thm:vmp.mono}
in Section \ref{sect:gen.p}.  The fact that
 $V_m(x)$ is decreasing in $m$, while $m V_m(x)$ is increasing
gives an indication of the delicate behavior of $V_m$.

 \item[d)] For $m > -1/2$, the definition of $V_m(x)$ can be
extended to $x = 0$ and
\begin{eqnarray}\label{vm0}
V_m(0) = \frac{\Gamma(m+\half)}{\Gamma(m+1)}.
\end{eqnarray}
For integer $m$, this becomes
\begin{eqnarray}\label{vm0.int}
V_m(0) = \frac{(2m)!}{2^{2m}(m!)^2} \sqrt{\pi} = \frac
  {1 \cdot 3 \cdot 5 \ldots (2m-1)}{2 \cdot 4 \cdot 6 \ldots (2m)} \sqrt{\pi}
\end{eqnarray}
while for large $m$ Stirling's formula implies
\begin{eqnarray}\label{vm0.asymp}
   V_m(0)  \approx \left( \frac{m-\half}{m} \right) ^m
     \left( \frac{e}{m} \right) ^{1/2} \approx \frac{1}{\sqrt{m}}
\end{eqnarray}
which is consistent with property (a).  Boyd \cite{By2,M} has
proved the more precise estimates
\begin{eqnarray}
   \frac{\sqrt{m + \frac{3}{4} + \frac{1}{32m + 48} }}{m + \frac{1}{2}}
  < V_m(0) < \frac{1}{\sqrt{m + \frac{1}{4} + \frac{1}{32m + 32} }}
\end{eqnarray}

  \item[e)] For all $m \geq 0$, $V_m$ satisfies the differential equation
\begin{eqnarray}\label{vmdiffeq}
  V_m^{\prime}(x) = 2x \left( V_m - V_{m-1} \right) .
\end{eqnarray}
This can easily be verified using integration by parts in (\ref{vmequiv2}).

 \item[f)] For each fixed $m \geq 0$, $V_m(x)$ is  decreasing in $x$.

This follows directly from  (b) and (e).

 \item[g)] For $a > 0$, the expression $aV_m(ax)$ increases with $a$.  Hence
 $aV_m(ax) > V(x)$ when $a > 1$ and $aV_m(ax) < V(x)$ when $a < 1$.

This property follows easily from the definition (\ref{vmdef}) or
(\ref{vmequiv2}) and the observation that
$\ds{ \frac{a}{\sqrt{a^2x^2 + u}} =  \frac{1}{\sqrt{x^2 + \frac{u}{a^2}}}}$
is increasing in $a$.
It is used in the proof of Theorem \ref{thm:vmp.mono}
and is important in the study of one-dimensional models for
 atoms in magnetic fields in which the electron-electron
interaction takes the form of convex combinations of
$\frac{1}{\sqrt{2}} V_m \left( \frac{ \left| x_j - x_k \right|}{\sqrt{2}}
\right).$

 \item[h)]  $V_0(x)$ is convex in $x > 0$; however, $V_m(x)$ is {\em not}
convex when $m > \half$.

For $m= 0$,  the differential equation
(\ref{vmdiffeq})  becomes
$V_0^{\prime}(x) = 2[xV_0 -1]$.  Since
$xV_0 = \int_0^{\infty} \frac{e^{-u}}{\sqrt{1+ u/x^2}} du$ is increasing
for $x > 0$, it follows that  $V_0(x)$ is convex.
 When $m > \half$ it follows from (\ref{vmdiffeq}) and (b) that
 $\lim_{x \rightarrow 0} V_m^{\prime}(x) = 0.$  Since
$V_m^{\prime}$ is negative,
$V_m^{\prime}$ must decrease on some small interval $(0,x_0).$
One can also show $\lim_{x \rightarrow \infty} V_m^{\prime}(x) = 0$,
so that one expects that there is an $x_1$ such that
$V_m$ is concave on $(0,x_1)$ and convex on $(x_1,\infty).$
    In Section \ref{sect:recurs} we will see that the convexity
is recovered for the averaged potential $V_m^{\av}$.

 \item[i)]  For integer $m$, $1/V_m(x)$ is convex in $x > 0.$

This will be proved in Section \ref{sect:rat.bnd} as Theorem
\ref{thm:convex.m}. For large
$x$,  $1/V_m(x) \approx x$ so that the deviation from linearity is very small
and the second derivative close to zero.  This makes the proof quite delicate
and lengthy.

The convexity of $1/V_m(x)$ can be rewritten as
\begin{eqnarray*}
   \frac{1}{ \half V_m \left( \frac{x+y}{2} \right) } \leq
      \frac{1}{V_m(x)} +  \frac{1}{V_m(y)} .
\end{eqnarray*}
Using property (g) with $a = \half$, one easily finds that the
 convexity of $1/V_m(x)$ implies
\begin{eqnarray*}
   \frac{1}{ V_m(x+y) } \leq
      \frac{1}{V_m(x)} +  \frac{1}{V_m(y)}.
\end{eqnarray*}
This subadditivity inequality plays the role of the triangle
inequality in applications.  (See, e.g. \cite{BR1}.)

  \item[j)]  Asymptotic estimates:

For large $x$,   it follows from property (a) that
\begin{eqnarray}
\frac{m}{2 (x^2+m)^{3/2}}
\leq  \frac{1}{x} - V_m(x) < \frac{m+1}{2x^3}
\end{eqnarray}
The asymptotic expansion
\begin{eqnarray}
 V_m(x) = \frac{1}{x} - \frac{m+1}{2x^3} + \frac{3(m+2)(m+1)}{8x^5}
      +  O \left(\frac{1}{x^7} \right)
\end{eqnarray}
can easily be obtained from (\ref{vmequiv2}).   For details
let $p = 2$ in the proof of Proposition \ref{prop:asymp.vmp}
which gives a similar expansion for $p >1$.

It follows from properties (a) and (c) that $V_m(x)$
decreases monotonically to zero for each fixed $x$ as $m \raw \infty$.
In fact, since $V_m(x)$ is decreasing in $x$ for all $m$, it suffices
to show this for $x = 0$ which is easy since $V_m(0) < m^{-1/2}$.

\item[k)]  The Fourier transform is given by
\begin{eqnarray}
\widehat{V}_m (\xi ) \equiv
\frac{1}{\sqrt{2\pi}} \int_{-\infty}^{\infty} V_m(x) \, e^{-ix\xi} dx =
\frac{4^{m+1}}{\sqrt{2 \pi}}\int _0 ^{\infty }
 \frac{ s^m \, e^{-s } }{(|\xi |^2 + 4s )^{m+1}}\,  ds
\end{eqnarray}
 This follows from  (\ref{vmdef})  and the standard formula
(e.g., see (v) on p. 131 of \cite{Stein})
$\ds{~{\cal F}\left( \frac{1}{\sqrt{|x|^2+|w|^2}} \right)(\xi) =
 \frac{1}{\sqrt{2 \pi}} \int_0^{\infty}
   \frac{1}{s} e^{- \half \left(s + |w|^2 |\xi|^2/s\right)} \,  ds~}$
after a change in the order of integration.

\end{itemize}

\subsection{Convexity Summary}
For large $x$, $V_m^p(x) \approx x^{1-p}$ which is convex in $x$
for $p > 1$ and concave for $p < 1$.  For $m > -\half$ these
convexity properties can not be extended to $V_m^p(x)$ on all of
$(0,\infty)$;  they would be inconsistent with the differential
equation and monotonicity properties in Proposition
\ref{prop:vmp.props}.  However, as discussed  after Proposition
\ref{prop:Vav.conv}, the averaged potentials $V_{\av}^{p,N}$ have
the same convexity as $x^{1-p}$ on the half-line.

The convexity of $1/V_m^p(x)$  is the motivation for Sections
\ref{sect4} and \ref{sect:rat.bnd.pfs}.  This question is already
delicate for $p = 2$ and its verification becomes increasingly
difficult for larger $p$.  Although, as discussed in Section
\ref{sect:genp}, we have evidence that convexity holds for all $p
\geq 2$,  our methods give this result only for a limited range
of $p$. Moreover, because our proof is inductive, we have
established convexity and ratio bounds of Section
\ref{sect:rat.bnd} only for integer $m$. It would be interesting
to find another approach which would extend these results to
non-integer $m$ and all $p \geq 2$.  Since $1/V_m^p(x) \approx
x^{p-1}$ which is concave for $1 < p < 2$, we can not expect
convexity of $1/V_m^p$ in this range.

As discussed above, one important consequence of the convexity of
$1/V_m(x)$ is an  analogue of the triangle inequality.  For all
values of $p$ we have $[V_m^p(x)]^{\frac{1}{1-p}} \approx x$ for
large $x$ which suggests a triangle inequality
of the form
\begin{eqnarray*}
 [V_m^p(x+y)]^{\frac{1}{1-p}} \leq
        [V_m^p(x)]^{\frac{1}{1-p}} + [V_m^p(y)]^{\frac{1}{1-p}}.
\end{eqnarray*}
It would be interesting to know the range of $p$ (and $m$) for
which this holds. For $p > 2$ the convexity of $1/V_m^p(x)$
implies only the weaker inequality
\begin{eqnarray*}
 [V_m^p(x+y)]^{\frac{1}{1-p}} \leq
        2^{\frac{p-2}{p-1}}\left([V_m^p(x)]^{\frac{1}{1-p}} +
[V_m^p(y)]^{\frac{1}{1-p}}\right).
\end{eqnarray*}

Finally, one could also ask if $V_m(x)$ is convex in $m$. In
particular, is $2V_m(x) \leq V_{m+1}(x) + V_{m-1}(x)$ or,
equivalently by (\ref{vmdiffeq}), is $V_m^{\prime}(x)$ increasing
in $m$?

\section{General $p$}\label{sect:gen.p}

We now study the basic properties of $V_m^p$ in detail.
As one would expect from $V_m^p(x) \approx x^{1-p}$, the behavior
of $V_m^p$ is often quite different for $p > 1 $ and $p < 1$.
At the boundary, $p =1$, $V_m^1(x) = 1$ for all $x$.
Proposition \ref{prop:vmp.in.p} describes the monotonicity and limiting
behavior of $V_m^p(x)$ as $p$ varies with $m$ and $x$ fixed.
Proposition \ref{prop:vmp=1/n} gives a simple expression
for $V_m^p$ in the special case that $1/p$ is an integer.

The next four results generalize properties of Section
\ref{subsect:propvm} to general $p$.  Proposition \ref{prop:vmp.ineq}
generalizes the inequalities from property (a);
Proposition \ref{prop:vmp.props} generalizes properties (d), (e), (f), and (g);
and, Theorem \ref{thm:vmp.mono} extends the monotonicity properties (b) and
(c).
Moreover, the proof of monotonicity for non-integer jumps is provided
here.  Finally, Proposition \ref{prop:asymp.vmp} gives the asymptotic
behavior of $V_m^p(x)$ for large $x$ when $p > 1$.

\begin{prop}\label{prop:vmp.in.p}
Let $ m > -1$ and $x > 0$ be fixed. Then
\begin{itemize}
\item[\rm{(i)}] $\lim_{p \rightarrow 0}V_m^p(x)=\infty$.
\item[\rm{(ii)}]
For all $x \geq 1$,
$V_m^p$ is decreasing in $p$.  Moreover,
\begin{eqnarray*}
 if ~x > 1,~~
\lim_{p \raw \infty}V_m^p(x) & = & 0, ~~~ and  \\
 if ~x = 1,~~  \lim_{p \raw \infty}
V_m^p(1) & = & \frac{1}{\Gamma(m+1)}\int_0^\infty
\frac{u^me^{-u}du}{1+u}.
\end{eqnarray*}
\item[\rm{(iii)}] For all  $0 < x < 1$ and $m > 0$,
$\lim_{p \rightarrow \infty}V_m^p(x)=\frac{1}{m}$.
\end{itemize}
\end{prop}

\pf We use the expression (\ref{vmp2}) for $V_m^p$.
\vskip 5mm
\par
(i) Since $\lim_{p \rightarrow 0} (x^p+u)^{1/p} = \infty$ for $x > 1$,
$$\lim_{p \rightarrow 0}V_m^p(x)\geq \lim_{p
\rightarrow 0} \frac{1}{\Gamma(m+1)}
\int_1^3 \frac{u^me^{-u}(x^p+u)^\frac{1}{p}du}{x^p + u}=\infty.$$
\par
(ii) Differentiating (\ref{vmp2}) yields
\begin{eqnarray}\label{vmp.diff.p}
\lefteqn{\frac{d}{dp} V_m^p(x) = \frac{1}{\Gamma(m+1)}\int_0^\infty
\frac{u^me^{-u}}{p^2(x^p+u)^{\frac{2p-1}{p}}} ~\times } \\
& ~ & ~~~~~~~~~~~~~~~~~~~  \nonumber
  \left[ \frac{}{}(1-p)x^p \ln(x^p)-(x^p+u)\ln(x^p+u) \right] du.
\end{eqnarray}
For $x=1$ or $p = 1$, the first term in square brackets above is zero
leaving a quantity which is clearly negative.
When both $x > 1$ and $p > 1$ both terms in (\ref{vmp.diff.p}) are
clearly negative.   When $x > 1$ and $0 < p < 1$,
the quantity in square brackets in (\ref{vmp.diff.p})
 is negative since
\begin{eqnarray*}
\lefteqn{(1-p) x^p\ln(x^p) - (x^p+u)\ln(x^p +u)} ~~~~~~~~~~~ \\
 & \leq & x^p\ln(x^p) - (x^p+u)\ln(x^p +u) < 0.
\end{eqnarray*}
The last inequality follows from the fact that the function
$f(w) = w \ln w $ is increasing for $w > 1$.
Thus,  $\frac{d}{dp}V_m^p(x) \leq 0$
for all $x \geq  1$ and for all $p \in (0,\infty)$.
\par
(iii) Since $\lim_{p \rightarrow \infty} (x^p+u)^{1/p} = 1$ for $x < 1$,
$$\lim_{p \rightarrow \infty}V_m^p(x)=
\frac{1}{\Gamma(m+1)}\int_0^\infty
u^{m-1}e^{-u}du=\frac{1}{m}.   \qed$$

\noindent{\bf Remark:} The behavior for $x < 1$ depends upon $m$.
For ``small" $m$,  there is a $p_0$ such that
$V_m^p$ decreases (below $1$) on $(0,p_0)$ and increases on $(p_0,
\infty)$ to
$\frac{1}{m}$. For ``big" $m$ , $V_m^p$ simply decreases  to $\frac{1}{m}$ as
$p$ increases.

The next result shows that  in the special case that $1/p$
is an integer $V_m^p$ reduces to
a polynomial in $x^p = x^{1/n}$ of degree $n-1$.
\begin{prop} \label{prop:vmp=1/n}
For $n \in {\bf N}$, $n \geq 2$,
$$ V_m^{1/n}(x)=\frac{1}{\Gamma (m+1)}\sum_{k=0}^{n-1}{n-1 \choose k}
 \Gamma (m+n-k)  ~x^{k/n}  $$
for all $x \geq 0$ and $m > -1$.
\end{prop}
\pf It follows from (\ref{vmp2}) that for $p = \frac{1}{n}$
\begin{eqnarray*}
 V_m^{1/n}(x)=\frac{1}{\Gamma (m+1)}
    \int_0^\infty \left(x^{1/n} + u \right)^{n-1} u^m e^{-u}  du.
\end{eqnarray*}
When $n$ is an integer $\geq 2$, the result then follows
easily from the binomial expansion applied to
$(x^\frac{1}{n} + u )^{n-1}$ and the definition of
the $\Gamma$-function.

\begin{prop}\label{prop:vmp.ineq}
For all $x >0$
\begin{eqnarray*}
\frac{1}{(x^p+m+1)^{\frac{p-1}{p}}} \leq V_m^p(x) \leq
\frac{1}{(x^p+m)^{\frac{p-1}{p}}}  ~~~{\rm for} ~ p > 1,
\end{eqnarray*}
where the first inequality holds for $m > -1$ and the second for $m \geq 0$.
\begin{eqnarray*}
(x^p+m+1)^{\frac{1-p}{p}} \geq V_m^p(x) \geq
(x^p+m)^{\frac{1-p}{p}} ~~~{\rm for} ~ \frac{1}{2} \leq p <1,
\end{eqnarray*}
where the first inequality holds for $m > -1$ and the second for $m \geq 0$.
\begin{eqnarray*}
V_m^p(x) \geq (x^p+m+1)^{\frac{1-p}{p}} ~~~{\rm for} ~0 < p \leq \frac{1}{2},
\end{eqnarray*}
and the inequality holds for $m > -1$.
\end{prop}

\pf
The proofs are done using Jensen's inequality as in property (a) of
Section \ref{subsect:propvm}.

\begin{prop}\label{prop:vmp.props} For all $x >0$
\begin{itemize}
\item[\rm{(i)}] For $m > - \frac{1}{p}$, $V_m^p(0)$ is defined and
$~\ds{V_m^p(0)= \frac{\Gamma(m+\frac{1}{p})}{\Gamma(m+1)}}.$
\item[\rm{(ii)}]For all $m \geq 0$, $x > 0$, $V_m^p$ satisfies the differential
equation
\begin{eqnarray}\label{diffeq.vmp}
\dffx V_m^p(x)=px^{p-1}\left(V_m^p(x)-V_{m-1}^p(x)\right).
\end{eqnarray}
\item[\rm{(iii)}]
For all $m > -1$, $V_m^p$ is decreasing in $x$, if $p >1$, identically
equal to $1$ for all $x$, if
$p=1$, increasing in $x$, if $p<1$.
\item[\rm{(iv)}] Let $m > -1$ and $x > 0$.
For $a > 0$, the expression $a^{p-1}V_m^p(ax)$ increases in $a$, if $p >1$
and decreases in $a$ if
$p<1$.
\end{itemize}
\end{prop}

\pf
The proofs are straightforward extensions of those given in Section
\ref{subsect:propvm}. In (iii), one can verify that $V_m^p$ is also
increasing  for $0 < p <
\frac{1}{2}$ by computing the derivative directly.

\begin{thm}\label{thm:vmp.mono}
For each fixed $x > 0$, and for $m$ in the region $m > -1$,
\begin{itemize}
\item[\rm{(i)}] $V_m^p(x)$ is strictly decreasing in $m$ for $p>1$ and
strictly increasing in $m$ for $p <1$.
\item[\rm{(ii)}] $mV_m^p(x)$ is strictly increasing in $m$ for $p >1$ and
strictly decreasing in $m$ for $p <1$.
\end{itemize}
\end{thm}

\pf To prove (i) we differentiate (\ref{vmp2}) to get
\begin{eqnarray}\label{eq:diffvm.wrtm}
  \frac{d~}{dm} V_m^p(x) =
  \frac{1}{\Gamma(m+1)}\int_0^\infty \frac{u^m \ln u~ e^{-u} }{(x^p +
u)^{1-\frac{1}{p}}} du
      - V_m^p(x) \frac{\Gamma^{\prime}(m+1)}{\Gamma(m+1)}
\end{eqnarray}
Using the same procedure as that used (see, e.g. \cite{AAR,L})
to obtain the standard integral representation
\begin{eqnarray}\label{eq:psirep}
 \psi(z) \equiv  \frac{ \Gamma^{\prime}(z)}{\Gamma(z)} =
 \int_0^{\infty} \left[ \frac{e^{-s}}{s} -\frac{1}{s(1+s)^z} \right] ds
\end{eqnarray}
one finds
\begin{eqnarray*}
\lefteqn{  \frac{1}{\Gamma(m+1)}
    \int_0^\infty \frac{u^m \ln u ~ e^{-u} }{(x^p + u)^{1-\frac{1}{p}}} du} \\
& =&\frac{1}{\Gamma(m+1)} \int_{s=0}^{\infty}   \frac{ds}{s}
    \int_0^{\infty} [e^{-s} - e^{-su}] \frac{e^{-u} u^m}{(x^p +
u)^{1-\frac{1}{p}}} du \\
   & = &V_m^p(x) \int_0^{\infty} \frac{e^{-s}}{s} ds -
           \frac{1}{\Gamma(m+1)}\int_0^{\infty}
\frac{ds}{s(s+1)^{m+\frac{1}{p}}}
      \int_0^{\infty}  \frac{e^{-w}w^mdw}{[x^p(s+1)+w]^{1-\frac{1}{p}}} \\
  & = & V_m^p(x) \int_0^{\infty} \frac{e^{-s}}{s} ds -
   \int_0^{\infty}  V_m^p(x(s+1)^\frac{1}{p})
\frac{ds}{s(s+1)^{m+\frac{1}{p}}},
\end{eqnarray*}
where we made the change of variable $w = (s+1) u$ to obtain
$V_m^p(x(s+1)^\frac{1}{p})$.
Now we use  Proposition \ref{prop:vmp.ineq} (iv) with $a = (s+1)^\frac{1}{p} >
1$ to obtain
\begin{eqnarray*}
\frac{1}{\Gamma(m+1)}
    \int_0^\infty \frac{u^m \ln u ~ e^{-u} }{(x^p + u)^{1-\frac{1}{p}}} du
  & \leq &
 \int_0^{\infty} V_m^p(x)
   \left(\frac{e^{-s}}{s}-\frac{1}{s(s+1)^{m+1}}\right) ds \\
  & = &   V_m^p(x) ~ \psi(m+1)
\end{eqnarray*}
when $p > 1$.  For $p < 1$, Proposition \ref{prop:vmp.ineq} (iv) gives
the inequality in the opposite direction.  Hence
inserting the result in (\ref{eq:diffvm.wrtm}) yields
\[\frac{d~}{dm} V_m(x) \left\{\begin{array}{ll}
< 0 & \mbox {if $ p >1$}\\
> 0 & \mbox {if $ p >1$}.\\
\end{array}
\right. \]
To prove (ii) it is slightly more convenient to consider
the logarithmic derivative $\frac{d~}{dm} \ln \left[ mV_m^p(x) \right]$
and show that it is  positive for $p>1$ and negative for $p < 1$.
Proceeding as above, we find for $p > 1$
\begin{eqnarray*}
\lefteqn{  \frac{d~}{dm} \ln \left[ mV_m^p(x) \right] =
\frac{1}{m} + \frac{\frac{d~}{dm} \left[ V_m^p(x) \right]}{V_m^p(x)}} ~~~~ \\
& = & \frac{1}{m} +  \int_0^{\infty} \frac{e^{-s}}{s} ds
    -  \frac{1}{V_m^p(x)} \int_0^{\infty}
    \frac{V_m^p \left[ x(s+1)^{1/p}\right]}{s~ (s+1)^{m+\frac{1}{p}}} ds
   - \psi(m+1) \\
& < &  \frac{1}{m} +  \int_0^{\infty} \frac{e^{-s}}{s} ds
     - \int_0^{\infty} \frac{1}{s (s+1)^m} ds
    - \psi(m+1) \\
& = &\frac{1}{m} + \psi(m) - \psi(m+1) \\
& = &\frac{1}{m} - \int_0^{\infty} \frac{1}{(s+1)^{m+1}} = 0
\end{eqnarray*}
where we have used (\ref{eq:psirep}) and
the following inequality with $a = (s+1)^{1/p}$.
\begin{eqnarray}
V_m^p(ax)  \left\{ \begin{array}{c} < \\ > \end{array} \right\}
 a V_m^p(x) ~~{\rm for}~~
\left\{ \begin{array}{c} p >1 \\ p <1 \end{array} \right\}.
\end{eqnarray}
for all $a \geq 1$.  This is easily verified
and implies that the inequality proved above for
$\frac{d~}{dm} \ln \left[ mV_m^p(x) \right]$ is reversed when $p < 1$.
\qed

The following result gives the asymptotic behavior of $V_m^p(x)$
for large $x$.
\begin{prop}\label{prop:asymp.vmp}  For $p > 1$, $V_m^p(x)$ has the
asymptotic expansion
\begin{eqnarray*}
\frac{1}{x^{p-1}} - \frac{(p-1)(m+1)}{p ~x^{2p-1}} +
\frac{(2p^2 -3p +1)(m^2 + 3m + 2)}{2p^2 ~x^{3p-1}}
      +  O \left(\frac{1}{x^{4p-1}} \right)
\end{eqnarray*}
\end{prop}
\pf This follows from (\ref{vmp2}) since
\begin{eqnarray*}
\lefteqn{V_m^p(x) = \frac{1}{\Gamma(m+1) ~ x^{p-1}} \int_0^\infty
          \frac{u^m e^{-u} }{ \left(1 + \frac{u}{x^p}
\right)^\frac{p-1}{p} } du }~~
\\
  & = &  \frac{1}{\Gamma(m+1)~x^{p-1}}\int_0^\infty  u^m
e^{-u}
  \left[ 1 - \frac{(p-1)~u}{p~x^p}  +
    \frac{(p-1)(2p-1)~u^2}{2p^2~x^{2p}}  +  \ldots
\right] du
\\
  & = &\frac{1}{x^{p-1}} \left[ 1 -
  \frac{(p-1)~\Gamma(m+2)}{p~\Gamma(m+1)~x^p}
+  \frac{(2p^2 -3p +1)~\Gamma(m+3)}{(2p^2)~\Gamma(m+1)~x^{2p}}
    +  O \left(\frac{1}{x^{3p}} \right) \right]
\end{eqnarray*}

\bigskip

\section{Recursion Relations and their Consequences}\label{sect:recurs}

\subsection{Recursion relations for $V_m^p$}\label{subs:recurs}

Although the case $p=2$ is of primary interest in applications, we
continue to study general $p$ is this section, as the proofs for
general $p$ are identical to those for $p=2$.  In these recursions,
our convention that $V_{-1}^p(x)=x^{1-p}$ plays an important role.

\begin{prop}\label{prop:recurs}
 For all $m \in {\bf R}$, $m \geq 1$, for all $x > 0$,
\begin{eqnarray}\label{vmprecurs}
  V_m^p(x) = \frac{1}{m}
   \left[ (m-1+ \frac{1}{p} -x^p)V_{m-1}^p(x)+ x^pV_{m-2}^p(x) \right].
\end{eqnarray}
\end{prop}
\pf
For $m=1$, one gets that
$$V_1^p(x)=pe^{x^p}\left[(\frac{1}{p}-x^p)\int_x^\infty
e^{-t^p}dt+x\right]=(\frac{1}{p}-x^p)V_0^p(x)+x^{p-1}V_{-1}^p(x).$$
For $m > 1$, using (\ref{vmp1}) and integration
by parts, we find
\begin{eqnarray*}
 V_m^p(x) & = & \frac{p e^{x^p}}{\Gamma(m+1)}\int_x^\infty  (t^p-x^p)^{m-1}
       (t^p-x^p)e^{-t^p} dt \\
& = & \frac{p e^{x^p}}{m\Gamma(m)} \left[
  (-x^p\int_x^\infty e^{-t^p} (t^p-x^p)^{m-1} dt  \frac{}{}  \right.  \\
      &~&~\left. \frac{}{}
  +  \frac{1}{p}\int_x^\infty
e^{-t^p}((t^p-x^p)^{m-1}+(m-1)pt^p(t^p-x^p)^{m-2} t^2)dt \right]\\
& ~ & \\
& = &  \frac{p e^{x^p}}{m\Gamma(m)} \left[(m-1+\frac{1}{p} -x^p)
  \int_x^\infty
e^{-t^p}(t^p-x^p)^{m-1} dt   \frac{}{}  \right.  \\
      &~&~~~~~~~~~~~~~~~~~~~~~~+ \left. \frac{}{} (m-1) x^p
   \int_x^\infty  (t^p-x^p)^{m-2} e^{-t^2} dt \right]\\
& = & \frac{1}{m} \left[ (m-1+\frac{1}{p} -x^p) V_{m-1}^p(x) +
   + x^p V_{m-2}^p(x) \right]
\end{eqnarray*}

 Repeated application of (\ref{vmprecurs}) gives a useful corollary.
For $m \in {\bf R}$, let $\fm$ denote the ``floor'' of
$m$, i.e., the largest natural number less than or equal to $m$.
\begin{cor} Let $m \in {\bf R}$, $m \geq 1$ and let  $n \in {\bf N}$ such
that  $n \leq \fm$. Then
\begin{eqnarray}\label{vmp.recurs2}
 V_m^p(x) & = & \frac{1}{pm}
   \left[ (1-px^p)V_{m-1}^p(x) + V_{m-2}^p(x)  + \ldots  \frac{}{}  \right.  \\
      &~&~~~~~~~\left. \frac{}{}  + \ldots [p(m-n)+1] V_{m-n}^p(x) + px^p
V_{m-n-1}^p(x) \right].
\nonumber
\end{eqnarray}
In particular, if $m$ is a positive integer, then
\begin{eqnarray}\label{vmpcor1}
 V_m^p(x) = \frac{1}{pm} \left[ (1 - px^p) V_{m-1}^p(x)  +
   \sum_{k=0}^{m-2} V_k^p(x) + px^p V_{-1}^p(x) \right].
\end{eqnarray}
\end{cor}

The expression (\ref{vmpcor1}) is well-defined for $x = 0$. Putting $x=0$ and
using  Proposition
\ref{prop:vmp.props} (i), we obtain the (presumably well-known) identity
\begin{eqnarray}
  \frac{\Gamma(m + \frac{1}{p})}{\Gamma(m+1)} = \frac{1}{pm} \sum_{k = 0}^{m-1}
       \frac{\Gamma(k + \frac{1}{p})}{\Gamma(k+1)}
\end{eqnarray}

\subsection{Averaged potentials}

These recursion relations are quite useful for studying
 the average  of the first $N$ of the $V_m$.  For $N$ a positive integer,
we extend (\ref{def.veff}) to
\begin{eqnarray}\label{def.veffp}
V_{\av}^{p,N}(x) = \frac{1}{N} \sum_{m=0}^{N-1} V_m^p(x)
\end{eqnarray}
Note that for $p=1$, $V_{\av}^{1,N}(x)=1$, for all $x \geq 0$.\\

The next result follows immediately from (\ref{vmpcor1}).
\begin{cor}$\ds{~~~ V_{\av}^{p,N}(x)  = p V_N^p(x) - \frac{px^p}{N}
   \left[ V_{-1}^p(x) - V_{N-1}^p(x) \right]}.$
\end{cor}
For the important case $p = 2$, this reduces to
\begin{eqnarray}
 V_{\av}^N(x)  = 2 V_N(x) - \frac{2x^2}{N}
   \left[ V_{-1}(x) - V_{N-1}(x) \right]
\end{eqnarray}

The function $V_0(|x|)$ is
convex on $(0,\infty)$ but has a cusp at $x = 0.$
However, as discussed in property (f), for higher $m$ both the convexity
and cusp are lost.  Thus, for higher $m$, the $V_m$ are somewhat smoother
than one might want for one-dimensional approximations to the Coulomb
potential.  The next result, although straightforward, is important because
it implies that the averaged potentials $V_{\av}^N(x)$ retain the cusp and
convexity properties of $V_0$ near the origin.

\begin{prop} \label{prop:Vav.conv}
The function
$\ds{V_{\av}^N(x)}$ is convex for all $x >0$ and
$$\ds{\lim_{x \raw 0+} \frac{d~}{dx} V_{\av}^N(x)} = - \frac{2}{N} .$$
\end{prop}
\pf  Using (\ref{def.veff})  and (\ref{vmdiffeq}) one finds
\begin{eqnarray*}
 \frac{d~}{dx} N V_{\av}^N(x) = \sum_{m=0}^{N-1} 2x[V_m(x) - V_{m-1}]
   = 2x[V_{N-1} - V_{-1}] = 2x V_{N-1} - 2.
\end{eqnarray*}
Therefore, to show that $V_{\av}^N$ is convex, we need to show that
$$x V_{N-1}=\frac{1}{\Gamma(N)} \int_0 ^\infty
\frac{u^{N-1}e^{-u}du}{(1+\frac{u}{x^2})^\frac{1}{2}}$$
is increasing. This holds as, for $x > 0$, the function $h_u(x) =
[1 + u/x^2]^{-1/2}$ is
increasing.

\bigskip
Similarly, one can show that for $p > 1$, $V_{\av}^{p,N}$ is convex
on $(0,\infty)$ and
$${\lim_{x \raw 0+} \frac{d~}{dx} V_{\av}^{p,N}(x)} = - \frac{p}{N}. $$
For $p <1$, the derviative becomes infinite at the origin; however,
 concavity of  $V_{\av}^{p,N}$  on $(0,\infty)$ still holds.

\subsection{Polynomials defined by recursion}\label{sect:poly}

We now observe that by repeatedly using (\ref{vmprecurs}) to eliminate
the $V_m^p$ with the largest value of $m$ from (\ref{vmp.recurs2}) allows
us to write
$V_m^p$ in terms of the two ``lowest'' functions (e.g., $V_0^p$ and
$V_{-1}^p$ in the case of integer $m$) and that the coefficients
in such expressions define two classes of polynomials related to
confluent hypergeometric functions.  We discuss the properties of
these polynomials in some detail. First, we
make the statement above explicit.

\begin{cor}  \label{cor:vmpoly}
For $m \geq 1$ there are polynomials $P_m^p(y)$ and $Q_m^p(y)$
of degree $\fm$ such that
\begin{eqnarray}\label{eq:vmp.poly}
  V_m^p(x)  & = &   P_m^p(x^p) V_{m-\fm}^p(x) + x^p Q_{m-1}^p(x^p)
       V_{m-\fm-1}^p(x).
\end{eqnarray}
\end{cor}
In the case of integer $m$  (\ref{eq:vmp.poly}) becomes
\begin{eqnarray}
    V_m^p(x)  & = &   \nonumber
  P_m^p(x^p) V_0^p(x) + x^{p} Q_{m-1}^p(x^p) V_{-1}^p(x) \\
  & = &   P_m^p(x^p) V_0^p(x) + xQ_{m-1}^p(x^p) \label{eq:vmppoly.int}
\end{eqnarray}
where the second expression follows from our convention
$V_{-1}^p(x) = x^{1-p}$.  We define
$P_m^p(y) = 1$ for $m \in [0,1)$ and $Q_m^p(y) = 0$
for $m \in [-1,0)$. Then (\ref{eq:vmp.poly}) holds trivially for
$m \in [0,1)$.
 \pf The desired polynomials are defined recursively.  First let
\begin{eqnarray}\label{defpoly1}
 P_m^p(y) & = & \frac{1}{m} \left( m-1+ \frac{1}{p} - y\right)
   ~~~ \mbox{for} ~~~ m \in [1,2) ,      ~~~ \mbox{and}  \\
Q_m^p(y) & = & \frac{1}{m+1}  ~~~ \mbox{for} ~~~ m \in [0,1) .
\end{eqnarray}
Then (\ref{eq:vmp.poly}) holds because it is equivalent to (\ref{vmprecurs})
for $m \in [1,2) $.  For $m \geq 2$ define $P_m^p(y)$   by
\begin{eqnarray}\label{eq:Pmp.recurs}
   P_m^p(y)=\frac{1}{m} \left[ \left( m-1 + \frac{1}{p} -y \right)
       P^p_{m-1}(y) + y~ P^p_{m-2}(y) \right]
\end{eqnarray}
and for $m \geq 1$ define $Q_m^p(y)$  by
\begin{eqnarray}\label{eq:Qmp.recurs}
Q_m^p(y)  & = &  \frac{1}{m+1} \left[ \left( m + \frac{1}{p} -y \right)
       Q^p_{m-1}(y) + y~ Q^p_{m-2}(y) \right].
\end{eqnarray}
It is straightforward to use induction to check that (\ref{vmprecurs})
yields (\ref{eq:vmp.poly}).  \qed

\bigskip
 We now restrict ourselves to
$m \in {\bf N}$ and study these polynomials in more detail.  The first few
polynomials are given in the following Table.
\begin{eqnarray*}\begin{array}{ccc}
    m  &  P_m^p  &  Q_m^p  \\ ~ \\
    0  &  1   &   1 \\
    1   &   \frac{1}{p} - y    &    \half(1 + \frac{1}{p} - y) \\
     2 &  ~~~\half\left[ \left(y - \frac{1}{p} \right)^2 +
  \frac{1}{p}\right]~~~  &
  \frac{1}{3} \left[ y + \half(1+ \fp -y)(2 + \fp -y) \right]
   \end{array}\end{eqnarray*}

\bigskip
 The following useful results, which hold for $m \geq 1$, are easily
checked by induction. $B(x,y)=\frac{\Gamma(x)\Gamma(y)}{\Gamma(x+y)}$
is the Beta function.
\begin{eqnarray*}
   P_m^p(0) & = &\frac{ \Gamma(m+ \fp)}{\Gamma(m+1)\Gamma(\fp)} =
\frac{1}{m~B(m,\fp)} \\
Q_m^p(0) & = & \frac{ \Gamma(m+1+ \fp)}{\Gamma(m+2)\Gamma(\fp)} =
\frac{1}{(m+1)~B(m+1,\fp)}
\end{eqnarray*}
\begin{eqnarray}
P_m^p(y) & = & \frac{1}{m} \left[\frac{1}{p} \sum_{j=0}^{m-1}P_j^p(y)
-y~P_{m-1}^p(y)  \right],  \label{eq:polysumP} ~~~\hbox{and} \\
Q_m^p(y) & = & \frac{1}{m+1}\left[\frac{1}{p}\sum_{j=0}^{m-1}Q_j^p(y)
        -y~ Q_{m-1}^p(y) + 1]\right].
\end{eqnarray}

We now obtain two expressions for $\dffx V_m^p(x)$.  First, observe that using
(\ref{eq:vmppoly.int}) in (\ref{diffeq.vmp}) yields
\begin{eqnarray*}
 \dffx V_m^p(x) =  p x^{p-1}
   \left( \left[ P_m^p(x^p) - P_{m-1}^p(x^p) \right] V_0^p(x) +
   x \left[Q_{m-1}^p(x^p) - Q_{m-2}^p(x^p) \right]
       \right).
\end{eqnarray*}
Differentiating  (\ref{eq:vmppoly.int}) yields after some
simplifications
\begin{eqnarray*}
\lefteqn{ \dffx V_m^p(x) = p x^{p-1}  \left[
   \left(P_m^p\right)^{\prime}(x^p) V_0^p(x) +
       P_m^p(x^p) [ V_0^p(x) -  V_{-1}^p(x) ] \right.} ~~~~~~ \\
  & ~ & ~~~~~~~~~~~~~~\left.  + ~ x^p
\left(Q_{m-1}^p\right)^{\prime}(x^p) V_{-1}^p(x)   +
    \fp Q_{m-1}^p(x^p) V_{-1}^p(x) \right]
\end{eqnarray*}
where $\left(P_m^p\right)^{\prime}(y)$ denotes $\dffy P_m^p(y)$.
Equating these expressions yields
\begin{eqnarray}
 \lefteqn{ - \left[\left(P_m^p\right)^{\prime}(x^p) +  P_{m-1}^p(x^p)
    \right] V_0^p(x) = } ~ \label{diffeq.pT} \\
& &  \left[ x^p \left(Q_{m-1}^p\right)^{\prime}(x^p) -  P_m^p(x^p)
+ (\fp - x^p) Q_{m-1}^p(x^p) +  x^p Q_{m-2}^p(x^p)  \right]
V_{-1}^p(x)
\nonumber  \end{eqnarray}
This provides motivation for the following
\begin{lemma}\label{lemma:diffpoly}
For $m \in {\bf N}$, $m \geq 1$,
\begin{eqnarray}
 \dffy P_m^{p}(y) & = & - P_{m-1}^p(y) , ~~~\hbox{and}   \label{eq:diffpoly} \\
  y~ \dffy Q^p_{m-1}(y) & = &
     P_m^p(y) - (m+1)  Q_m^p(y) +  m Q^p_{m-1}(y)  \label{diffeq.Qmp}
\end{eqnarray}
\end{lemma}

\pf We first prove (\ref{eq:diffpoly}) by induction. It can be verified for
$m =1,2$ using the Table above.   Then using (\ref{eq:Pmp.recurs}), we find
\begin{eqnarray*}
\lefteqn{m \dffy P_m^{p}(y)} \\
& = &\left( m-1 + \ts{\frac{1}{p}} -y \right)
       \dffy P^p_{m-1}(y) - P^p_{m-1}(y) + y  \dffy P^p_{m-2}(y)  +
P^p_{m-2}(y) \\
  & = &- P^p_{m-1}(y) - \left( m-2 + \ts{\frac{1}{p}} -y \right)
P^p_{m-2}(y) -  y~P^p_{m-3}(y) \\ & = &  - m P_{m-1}^{p}(y).
\end{eqnarray*}
This implies that the coefficient of $V_0^p$ in (\ref{diffeq.pT}) is
identically zero. Therefore the coefficient of $V_{-1}^p$ must also
be  identically zero. Substituting $y = x^p$ gives
\begin{eqnarray*}
  y~ \dffy Q^p_{m-1}(y) & = & P_m^p(y) +
 \left(y  - \frac{1}{p} \right)Q^p_{m-1}(y)  - y Q^p_{m-2}(y)  \\
  & = & P_m^p(y) - (m+1) Q_m^p(y) +  m Q^p_{m-1}(y)
\end{eqnarray*}
where we used (\ref{eq:Qmp.recurs}).  ~~~\qed

Note that since the left side of (\ref{diffeq.Qmp}) is a polynomial of
degree $m-1$, this implies the coefficients of the $y^m$ terms in
$P_m^p$ and $(m+1) Q_m^p(y)$ are identical.   In fact, one
can use (\ref{eq:Pmp.recurs}) and (\ref{eq:Qmp.recurs}) to see
that the leading terms of $P_m^p$ is $(-1)^m y^m/m!$
and that of $Q_m^p$ is $(-1)^m y^m/(m+1)!$

A set of polynomials $\{ p_n(x) \}$belongs to the class known as Appell
polynomials \cite{BB} if they satisfy  $\dffx p_n(x) = p_{n-1}(x)$.
Therefore,  (\ref{eq:diffpoly}) implies that for each fixed $p$,
the set $(-)^m P_m^p(y)$ forms a family of Appell polynomials.

\bigskip

One can use (\ref{eq:diffpoly})  in (\ref{eq:Pmp.recurs}) to
replace $P_{m-1}^p$ and
$P_{m-2}^p$ by derivatives of $P_m^p$ and obtain a
second order differential equation satisfied by $P_m^p$.
This allows us to obtain
a relationship between the polynomials
$P_m^p$ and confluent hypergeometric functions, which we
denote  $~_1F_1(\alpha,\gamma,y)$.

\begin{thm} For $m \in {\bf N}$, $m \geq 1$, $P_m^p(y)$ satisfies the
differential equation
\begin{eqnarray}\label{DE:polyP}
 y \phi^{\prime\prime} (y) -
  \left( m-1 + \frac{1}{p} - y \right) \phi^{\prime}(y) -
  m\phi(y) =0.
\end{eqnarray}
\end{thm}

Standard techniques show that (\ref{DE:polyP}) has a polynomial solution
of the form $\phi(y) = \sum_{k=0}^m b_k y^k$ with
$b_k = - \frac{m+1 -k}{k(m+ \frac{1}{p} - k)} b_{k-1}$, $k \geq 1$, $b_0 \neq
0$ arbitrary, and a second solution of the form
$\phi(y) = \sum_{k=0}^{\infty} c_k y^{k+m + 1/p}$ with
$c_k =  - \frac{k -1 +1/p}{k(m+ \frac{1}{p} + k)}c_{k-1}$, $k \geq 1$,
$c_0 \neq 0$ arbitrary.
Since $P_m^p(0)  = \frac{1}{m~B(m,\fp)}$, we conclude that
$b_k = \frac{(-1)^k}{k!~(m-k)B(m-k,\fp)}$ and
$$P_m^p(y) =  \sum_{k=0}^m \frac{(-1)^k \Gamma(m + \fp -k)}
  {\Gamma(k+1)\Gamma(m + 1 -k)\Gamma(\fp)}~ y^k.$$

The restriction that $1/p$ be non-integer in the next result
is neither serious, nor unexpected, in view of Proposition \ref{prop:vmp=1/n}.
\begin{cor}\label{cor:conf.hyper} Let $ p \neq \frac{1}{n}$
for $n \in {\bf N}$. Then
\begin{eqnarray}\label{Pm:conf.hyper}
     P_m^p(y) = \frac{1}{m~B(m,\fp)}~ e^{-y}
  ~_1F_1\left(1 - \fp,1 - \fp-m,y \right)
\end{eqnarray}
\end{cor}
\pf Write $\phi(y) = e^{-y} \widehat{\phi}(y)$.  Then it follows from
(\ref{DE:polyP}) that
$\widehat{\phi}$ satisfies
\begin{eqnarray}\label{conf.hyper}
 y \widehat{\phi}^{\prime\prime} (y) -
  \left( m-1 + \frac{1}{p} + y \right) \widehat{\phi}^{\prime}(y) -
   \left(1 - \frac{1}{p}\right)\widehat{\phi} (y) = 0.
\end{eqnarray}
which has the form of the differential equation satisfied by the confluent
hypergeometric function.   Comparing the behavior of $P_m^p(y)$ near
$y=0$ with that of the well-known solutions to (\ref{conf.hyper})
suffices to complete the proof.
 \qed

\bigskip

 It is well-known \cite{Bate} that for real $\alpha$
and $\gamma$,  $~_1F_1(\alpha,\gamma,y)$ has at most
finitely many zeros on the real line.  Hence, the same
holds for $P_m^p$.  In fact, we can show that $P_m^p$
has no zeros when $m$ is even and exactly one when $m$ is odd.

To show this, it is convenient to introduce the new
variable $z =  \frac{1}{p}-y$ and write
$P_m^p(y) = \tilde{P}_m^p \left(\frac{1}{p}-y \right)$.  The first
few of these polynomials are $\tilde{P}_0^p(z) = 1$,
$\tilde{P}_1^p(z) = z$, and
$\tilde{P}_2^p(z) = \half\left(z^2 + \frac{1}{p}  \right)$.
\begin{lemma}\label{lemma:polyz} For $m \in {\bf N}$, $m \geq 2$,
the polynomials  $\tilde{P}_m^p(z)$ satisfy
\begin{itemize}
\item[\rm{(i)}] $\ds{\tilde{P}_m^p(z)=
  \frac{1}{m} \left[ (m-1+z)\tilde{P}_{m-1}^p(z)
   +\left(\frac{1}{p}-z \right)\tilde{P}_{m-2}^p(z) \right]}$,
\item[\rm{(ii)}]$\ds{\tilde{P}_m^p(z)=\frac{1}{m} \left[
   \frac{1}{p}\sum_{j=0}^{m-2}\tilde{P}_j^p(z) +
       z\tilde{P}_{m-1}^p(z)\right]}$, and
\item[\rm{(iii)}] $\frac{d}{dz}\tilde{P}_m^p(z)=\tilde{P}_{m-1}^p(z)$.
\end{itemize}
\end{lemma}
\pf Follows immediately from substitution in (\ref{eq:Pmp.recurs}),
(\ref{eq:polysumP}) and (\ref{eq:diffpoly}).

\begin{cor}\label{cor:pos.coef}
For $m \in {\bf N}$, $m \geq 0$, all coefficients
in the polynomials $\tilde{P}_m^p(z)$ are positive.
\end{cor}

\pf  This follows immediately from the explicit expressions
above for $\tilde{P}_m^p$ when $m = 0,1$ and part (ii) of
Lemma \ref{lemma:polyz}.

\begin{prop}\label{prop:roots} Let $m \in {\bf N}$, $m \geq 1$
\begin{itemize}
\item[~] If $m$ is even, $\tilde{P}_m^p(z) \geq 0$.
\par
\item[~]If $m$ is odd,  $\tilde{P}_m^p(z)$ has
exactly one root $z_m$.
\end{itemize}
Moreover, the roots form a strictly decreasing sequence  with
$-m+1 \leq z_m \leq 0$.
\end{prop}

\pf
First note that Corollary \ref{cor:pos.coef} implies that
$\tilde{P}_m^p(z) \geq 0$ for all $z \geq 0$ and $\lim_{z
\rightarrow \infty}\tilde{P}_m^p(z)=\infty$.
Now we claim that
\[\lim_{z \rightarrow - \infty}\tilde{P}_m^p(z)=\left\{\begin{array}{rl}
- \infty & \mbox{if m is odd}\\
\infty & \mbox{if m is even}
\end{array}
\right.\]
This can easily be verified by induction using part(i) of
 Lemma \ref{lemma:polyz} (i) above.
\par
For $m$ odd, $\tilde{P}_m^p(z)$ has at least one root $z_m$.
We now prove by induction that $z_m$ is the only root if $m$ is odd
and that $\tilde{P}_m^p(z) \geq 0$ for all $z$ if $m$ is even.
The induction hypothesis is easily seen to hold for $m=0,1$. Suppose
it is true up to $m-1$ and
consider $\tilde{P}_m^p$, with $m$ odd. Then,
by Lemma \ref{lemma:polyz} (iii),
$\frac{d}{dz}\tilde{P}_m^p(z)=\tilde{P}_{m-1}^p(z)$. Since $m-1$ is even,
by the induction hypothesis $\tilde{P}_{m-1}^p(z) \geq 0$.
 Thus $\tilde{P}_m^p(z)$ is increasing for all $z \in {\bf R}$, which
implies that for $m$ odd $\tilde{P}_m^p(z)$ has only one root.
Since  $\tilde{P}_m^p(z) > 0$ when $z > 0$ that one root must
satisfy $z_m \leq 0$.
\par
If $m$ is even, then $\frac{d}{dz}\tilde{P}_m^p(z)=\tilde{P}_{m-1}^p(z)$
which, by the  induction hypothesis,  has
exactly one root $z_{m-1} \leq 0$. Therefore $\tilde{P}_m^p(z)$ has a local
extremum at $z_{m-1}$.  Since
$\frac{d^2}{(dz)^2}\tilde{P}_m^p(z)=\tilde{P}_{m-2}^p(z)$ and $m-2$ is
even,   $\tilde{P}_{m-2}^p(z)\geq 0$ by the
induction hypothesis and $z_{m-1}$ is a local minimum
for $\tilde{P}_m^p(z)$.
By Lemma \ref{lemma:polyz} (i) we have
$$ m\tilde{P}_m^p(z_{m-1}) =
  \left(\ts{\ts{\frac{1}{p}}}-z_{m-1} \right) \tilde{P}_{m-2}^p(z_{m-1}) $$
since $\tilde{P}_{m-1}^p(z_{m-1})=0$. But by the induction hypothesis
$z_{m-1}\leq 0$ and $\tilde{P}_{m-2}^p(z_{m-1}) > 0$
so that  $\tilde{P}_m^p(z_{m-1}) > 0$ as required.

It remains to be shown that the roots are decreasing and bounded
below by $ -(m-1)$.   Both can be easily checked for   $m=1,3$
and then proved by induction using  Lemma \ref{lemma:polyz} (i).
We now let $m$ be odd.  Since $\tilde{P}_m^p(z)$
is increasing to show that $z_m > -m + 1$, it suffices to
show that  $\tilde{P}_m^p( -m +1) < 0$. For
$z = -m+1$ the recursion relation reduces to
$$ m\tilde{P}_m^p(-m+1) = \left(\fp +m-1\right)\tilde{P}_{m-2}^p(-m+1)$$
which is negative by the induction assumption
hypothesis that $z_{m-2} \geq -m+3$.
  To show that $z_{m} < z_{m-2}$ it
suffices to show that $\tilde{P}_m^p(z_{m-2}) > 0$.
But
$$m\tilde{P}_m^p(z_{m-2})= (m-1+z_{m-2}) P_{m-1}^p(z_{m-2})
 \geq 0$$
since $P_{m-2}^p(z_{m-2}) = 0$,  $z_{m-2} > -m + 3 > -m+1$, and
 $P_{m-1}^p(z)$ is positive.  \qed

We now restate the results above in terms of the behavior of the
original polynomials $P_m^p(y)$.
\begin{cor}\label{rootsy}Let $m \in {\bf N}$, $m \geq 1$.
Then
\begin{itemize}
\item[{\rm (i)}]
  If $m$ is even, then $P_m^p(y) \geq 0$ for all $y \in {\bf R}$.
\item[{\rm (ii)}] If $m$ is odd, then $P_m^p(y)$ has
exactly one root $y_m \geq 0$.
\item[ {\rm (iii)}] For all $m \geq 1$,
\begin{eqnarray*}\lim_{y \rightarrow - \infty}P_m^p(y) & = & \infty
\\ \lim_{y \rightarrow \infty}P_m^p(y) & = & \left\{\begin{array}{rl}
- \infty & \mbox{if m is odd}\\
\infty & \mbox{if m is even}
\end{array}
\right.
\end{eqnarray*}
\end{itemize}
\end{cor}

Although we were able to obtain an explicit expression for the
polynomials $P_m^p(x)$  relating them to confluent hypergeometric
functions and analyze their behavior in some detail, we do not
much information about $Q_m^p(x)$.  This is, at least in part,
because (\ref{diffeq.Qmp}) mixes $Q_m^p(x)$ and $P_m^p(x)$ and
does not lead directly to a differential equation for $Q_m^p(x)$.
It would be interesting to know more about the polynomials
$Q_m^p(x)$.

\section{Inequalities and Convexity} \label{sect4}

\subsection{Inequalities for $V_0(x)$}\label{sect:V0}

We first illustrate our strategy by proving a special class of
inequalities for $V_0$. The convexity of
$1/V_0(x)$ follows directly from the optimal upper bound
in this class as given in Theorem \ref{ineqv0} below.
Although, as discussed at the end of this section,  these inequalities
generalize to $V_m$ the resulting upper bound is not sufficient to
establish the convexity of $1/V_m$.  For this we need a bound on
the ration $V_m(x)/V_{m-1}(x)$.  Nevertheless, these simple
inequalities for $V_0$, which can also be interpreted as ratio bounds,
are of some interest in their own right in a variety of applications.
Because the geometric strategy  is also used in our
more complex proofs of ratio bounds, we think there is some merit in
presenting it first here.

We now define
\begin{eqnarray}\label{gkdef}
         g_k(x) = \frac{k}{(k-1)x + \sqrt{x^2+k}}.
\end{eqnarray}

\begin{thm}\label{ineqv0}
For $x \geq 0$
        \begin{eqnarray}\label{opt.ineq}
          g_{\pi}(x) \leq  V_0(x) < g_4(x)
        \end{eqnarray}
and these inequalities
         are optimal for functions of the form (\ref{gkdef})
        with equality only at $g_{\pi}(0) = V_0(0) = \sqrt{\pi}.$
\end{thm}
\pf
It is easy to see that the family of functions
        $g_k(x)$ is increasing in $k$ and that $ 0 < g_k(x) < 1/x.$
In order to prove that the upper bound is optimal, we first observe
        that  $ g_k^{\prime}(x) = -k [g_k(x)]^2
      [x + (k-1)\sqrt{x^2+k}~]/{\sqrt{x^2+k}}$ and
     $x g_k(x) -1 = -k g_k(x)/[x + \sqrt{x^2+k}~]$.
Then one can verify that
        \begin{eqnarray}\label{ineq:deriv}
\lefteqn{ g_k'(x) > 2[x g_k(x)  - 1] } ~~~~~~~~ \nonumber \\
        & \iff &  \frac{k}{\sqrt{x^2+k}}~
    \frac{(k-1)\sqrt{x^2+k} + x}{(k-1)x + \sqrt{x^2+k}} <
      \frac{2k}{x + \sqrt{x^2+k}} \nonumber \\
           & \iff & (k-2) x^2 + k (k-3) < (k-2) x \sqrt{x^2+k} \nonumber
        \\ & \iff & x^2 (k-2)(k-4) + k(k-3)^2 < 0
        \end{eqnarray}
        when $k > 3.$  We now restrict attention to $3 \leq k \leq 4$ and
        let $h_k(x) = g_k(x) - V_0(x).$  For $k=4$, the expression
        (\ref{ineq:deriv}) implies  $g_4'(x) < 2[xg_k(x)  - 1]$
        so that $h_4^{\prime}(x) < 2x h_4(x)$ for all $ x \geq 0$; whereas for
        $k < 4$ this holds only for
        $x < a_k = \sqrt{\frac{k(k-3)^2}{(k-2)(4-k)}}.$
        Since both $V_0(x)$ and $g_k(x)$ are positive and bounded above
        by $1/x$, their difference also satisfies
        $|h_k(x)| < 1/x \rightarrow 0.$

        For $k=4$, if $h_4(x) \leq 0$ for some $x > 0 $, then  $h_4'(x)< 2x
h_4(x)$
    is negative
        and thus $h_4$ is negative and strictly decreasing from a certain
$x$ on,
        which contradicts
        $\lim_{x \rightarrow \infty} h_4(x) = 0.$  Thus $h_4(x) > 0$ so
        that $g_4(x) > V_0(x)$, for all $x$.
        Now suppose that
        for some $k < 4$, $g_k$ is an upper bound, i.e. $h_k(x) \geq 0$
        for all $x \geq 0$.
        In particular, $h_k(x) \geq 0$ for all $x > a_k$.
        For $k < 4$, we find however
        that  $h_k^{\prime}(x) > 2 x h_k(x)$ holds for $x > a_k$.
 Thus we get $h_k(x) \geq 0$ and strictly increasing for all $x > a_k$
which
contradicts $\lim_{x \rightarrow \infty} h_k(x) = 0.$
        Thus the upper bound  can {\em not} hold when
         $x > a_k$ and $k < 4$. The lower bound also fails
         for $k > \pi$ since then
         $h_k(0) = g_k(0) - V_0(0) = \sqrt{k} - \sqrt{\pi} > 0.$

    To establish the improved lower bound $g_{\pi} \leq V_0(x)$
    we note that the argument above implies that
    $h_k(x)$ is negative for $x > a_k$ and $3 < k \leq \pi.$
    However for $k < \pi$ we have $h_k(0) < 0$ so that
    $h_k(x)$ is also negative for very small $x$.  If $h_k(x)$
     is ever non-negative, we can let $b$ denote the first
     place $h_k(x)$ touches or crosses the x-axis, i.e., $h_k(b) = 0$
     and $h_k(x) < 0$ for $x < b$.   Then $h_k$ must be increasing
     on some interval of the form $(x_0,b)$.  However, by the
     remarks above, $h_k(b) = 0$ implies $b \leq a_k$ so that
     $h_k^{\prime}(x) < 2 x h_k(x) < 0$ on $(x_0,b).$  Since this
     contradicts $h_k$ increasing on $(x_0,b)$, we must have
     $h_k(x) < 0$ for all $x \geq 0$ if $k < \pi.$

    Thus we have proved the lower bound
    $g_k(x) < V_0(x)$ on $[0,\infty)$  for $k < \pi.$
    Since $g_k$ is  continuous and increasing in $k$, it follows that
    $g_{\pi}(x) \leq V_0(x).$  To show that this inequality is strict
    except at $x=0$, note that the right derivative of $h_k$ at $0$
    satisfies $h_k'(0) = 2 - k$ so that  $h_{\pi}'(0) < 0$ and
    $h_{\pi}(x)$ is negative at least on some small interval $(0,x_1).$
    Then we can repeat the argument above to show that $h_{\pi}(x) < 0$
    if $x > 0.$  \qed

\bigskip

As discussed in \cite{BR1,SzW, W} the upper bound implies the convexity of
$1/V_0(x)$ on $(0,\infty)$;  in fact, it is not hard to
use the fact that  (\ref{vmdiffeq}) reduces to $\dffx V_0(x) = 2(xV_0 - 1)$
to see that the upper bound is equivalent to convexity.
It  was established independently by
        Wirth \cite{W}  and by Szarek and Werner \cite{SzW}.
(The latter actually proved slightly more by using (\ref{vmequiv1}) to
define an asymmetric extension of $V_0(x)$ to negative $x$.  They showed in
\cite{SzW} that this extension is convex for  $x >-\frac{1}{\sqrt 2}$.)

Both bounds in (\ref{opt.ineq}) are  sharper than
        the  inequalities of Komatsu \cite{IM,M}.
        The weaker lower bound $g_3(x) < V_0(x)$ was used in \cite{BR1} to show
        that the function $ \left[ 1/V_0(x) - x \right]^2/V_0(x)$ is decreasing
        for $x \geq 0$.
The lower bound $ g_{\pi}(x) \leq  V_0(x)$ was established earlier
by Boyd \cite{By1} as the optimal bound in a different class
of inequalities.  There is an extensive literature (see. e.g. \cite{M}) on
bounds for $V_0(x)$; however, the class of inequalities obtained
using functions of the form $g_k(x)$ does not seem to have been considered
before so that the optimality of bounds of this type for $k = \pi$ and
$k = 4$ seems new.

\bigskip

Mascioni  \cite{Mas} generalized the upper bound to $p \geq 2$
for which he showed
$$V_0^p(x) <
\frac{4p}{3px^{p-1}+\sqrt{p^2x^{2p-2}+8p(p-1)x^{p-2}}}
$$
and also showed that this implies convexity of $1/V_0^p(x)$
for $p \geq 2$.

\bigskip

In view of Property (a) of Section \ref {subsect:propvm},
it would seem natural
to try to generalize (\ref{ineqv0}) using functions of the form
\begin{eqnarray}\label{gkm.def}
         g_k^m(x) = \frac{k}{(k-1)x + \sqrt{x^2+m+k}}.
\end{eqnarray}
Note that the functions $g_k^m$ are increasing in $k$ and that
$\lim_{k \rightarrow \infty} g_{k}^m(x)=\frac{1}{x}$.
Therefore Property (a) implies that
\begin{eqnarray*}
g_{1}^m(x) \leq  V_m(x) <  \lim_{k \rightarrow \infty} g_{k}^m(x).
\end{eqnarray*}
As $g_k^m$ is continuous in $k$, there must exist $i_m$ and $j_m$ such that
\begin{eqnarray}\label{ineq.gkm}
g_{i_m}^m(x) \leq  V_m(x) < g_{j_m}^m(x).
\end{eqnarray}
However, we have not obtained explicit expressions for
$i_m$ and $j_m$.  One might expect that the optimal lower bound
occurs when  $i_m$ is chosen to satisfy $g_{i_m}^m(0) =  V_m(0)$.
However, numerical evidence shows that this is false; in fact, this
choice for $i_m$ does not even yield an inequality.

\subsection{Ratio Bounds}\label{sect:rat.bnd}

One of our main goals  is to show that the function $\frac{1}{V_m(x)}$ is
convex for integer $m \geq 1$.  The key to this is the realization that
(\ref{opt.ineq}) can also be rewritten to give bounds on the ratio
$V_0(x)/V_{-1}(x) = xV_0(x)$.  We now let
\begin{eqnarray}\label{def:Gkm}
    G_k^m(y) = \frac{ky}{ (k-1)y - m + \sqrt{(y+m)^2 + ky} }
\end{eqnarray}
and note that $ x g_k(x) = G_k^0(x^2)$ so that (\ref{opt.ineq}) is
equivalent to
\begin{eqnarray*}
  G_{\pi}^0(x^2) \leq x V_0(x) = \frac{V_0(x)}{V_{-1}(x)} < G_4^0(x^2).
\end{eqnarray*}
For integer $m > 0$, convexity of $\frac{1}{V_m(x)}$ can be
shown to be equivalent to
\begin{eqnarray*}
  R_m(x) \equiv \frac{V_m(x)}{V_{m-1}(x)} < G_4^m(x^2).
\end{eqnarray*}

\bigskip

In addition to this upper bound, we can show the
following
\begin{thm} \label{thm:optm.ineq} Let $m \in {\bf N}$, $m \geq 0$. Then the
inequalities
\begin{eqnarray}\label{eq:rat.bnd}
G_8^{m-1}(x^2) <   R_m(x)  < G_4^m(x^2)
\end{eqnarray}
hold and are optimal in $k$ for the class of functions of the form
$G_k^m(x^2)$.
\end{thm}

The upper bound is optimal in $k$ for all $m$.  The lower bound
is optimal in the sense that $8$ is the largest integer for
which the lower bound in (\ref{eq:rat.bnd}) holds for all $m$.
However, as we discuss at the end of Section  \ref{sect:opt}, for
fixed $m$ one can find $k(m)$ such that
$G_{k(m)}^{m-1}(x^2) <   R_m(x) $ holds with $k(m) > 8 $.

Since $G_k^m(y)$ is increasing in both $m$ and $k$, its behavior at
zero and infinity allows us to also draw some conclusions about the
optimality in $m$ of (\ref{eq:rat.bnd}).
 $R_m(0)=1-\frac{1}{2m}$ and $G_k^m(0)=1-\frac{1}{1+2m}$ for all $k$.
Therefore, $G_k^{\nu}(0) <   R_m(0)  < G_{k'}^{\mu}(0)$ implies
$\nu \leq m-\half$ and $\mu \geq m - \half$ for all  $k, k'$.
Thus, if we insist that $m$ be integer, there is no choice of $k$
which allows $m-1$ to be replaced by $m$ in the lower bound when
 $m > 0$, or $m$ by $m-1$ in the upper bound.  This argument does not,
however, rule out the possibility of bounds of the form
$G_k^{m-\half}(x^2) <   R_m(x)  < G_{k'}^{m-\half}(x^2)$.

To examine the behavior at infinity, note that
\begin{eqnarray*}
  G_k^m(y) & = & 1 - \frac{1}{2y} + \frac{4m + k + 2}{8y^2} +
    O\left( \frac{1}{y^3} \right), ~~~\hbox{and} \\
  R_m(\sqrt{y})  & = & 1 - \frac{1}{2y} + \frac{4m + 6}{8y^2} +
    O\left( \frac{1}{y^3} \right)
\end{eqnarray*}
where the asymptotic expansion for $R_m$ follows from Proposition
\ref{prop:asymp.vmp}.  It then follows that
$ R_m(\sqrt{y})  < G_{k'}^{\mu}(y)$ implies $\mu > m + 1 - \frac{k}{4}.$
Thus $m$ is optimal for the upper bound if $k \leq 4$ and any
attempt to decrease $m$ would require an increase in $k$.
Furthermore, $\mu = m$ implies $k \geq 4$ so that the
upper bound in (\ref{eq:rat.bnd})  is optimal in $k$.

We postpone the proof of Theorem \ref{thm:optm.ineq}, which requires
a lengthy computation even for the case $p = 2$, to the next
section.   Our proof uses induction on $m$.  Therefore, we are
able to establish (\ref{eq:rat.bnd}) and the theorems in the
next section
 only for  $m$ a positive integer.  We believe that they are also
true for non-integer $m.$  However, a proof would require either a
different method or independent verification of the
upper bound
for an initial range, such as $-1 < m < 0$.

\bigskip

The ratio $ R_m(x) $ is
of interest in its own right, and our results are sufficient to
establish that it is increasing in $x$ on $(0,\infty)$.
This is proved in the next section after
Theorem \ref{thm:convex.m}, which uses a similar argument.
\begin{thm} \label{thm:ratio.inc}
For $m \in {\bf N}$,
 the ratio $R_{m+1}(x)=\frac{V_{m+1}(x)}{V_m(x)}$ is
increasing in $x$.
\end{thm}

\subsection{Convexity of $1/V_m$}

We now prove some important
 consequences of Theorem \ref{thm:optm.ineq}.  The first
is
\begin{thm} \label{thm:convex.m}  For all $m \in {\bf N}$,
the function $1/V_m(x)$ is convex on $[0,\infty)$.
\end{thm}
\pf We need to show that
\begin{eqnarray}\label{eq:cnvx.cond}
\left(\frac{1}{V_m(x)}\right)^{\prime\prime}=
\frac{2[V_m(x)^{\prime}]^2- V_m(x)(V_m(x))^{\prime\prime}}{V_m(x)^3} > 0.
\end{eqnarray}
It follows from the differential equation (\ref{vmdiffeq})
and the recursion relation (\ref{vmprecurs}) that
\begin{eqnarray*}
V_m(x)^{\prime\prime}=2 \left[ V_m(x)(1+2m+2x^2)-2V_{m-1}(x)(x^2 + m) \right]
\end{eqnarray*}
so that (\ref{eq:cnvx.cond}) holds if and only if
$$[V_m(x)]^2(1+2m-2x^2)+
2V_{m-1}(x)V_m(x)(3x^2-m)-4x^2[V_{m-1}(x)]^2 \leq 0.$$
After division by $[V_{m-1}(x)]^2$ this can be rewritten as
$P[R_m(x)] \leq 0$ where
\begin{eqnarray*}
P(z)=z^2(1+2m-2x^2)+2z(3x^2-m)-4x^2 .
\end{eqnarray*}
Writing the roots of $P(z) = Az + 2Bz + C$ in the non-standard
form $\frac{-C}{B \pm \sqrt{B^2 - AC}}$, we find that $G_4^m(x^2)$
is either the smaller of two positive roots (when $x^2 >m+\half$) or
the only positive root (when $x^2 < m+\half$).  Since $P(0) < 0$,
in both cases, we can conclude that
$$z < G_4^m(x^2) \hspace{3mm}\mbox{implies} \hspace{3mm} P(z) < 0.$$
Therefore, it follows from the upper bound in Theorem \ref{thm:optm.ineq}
that $P[R_m(x)] < 0$, hence (\ref{eq:cnvx.cond}) holds. \qed
\bigskip

\noindent{\bf Proof of Theorem \ref{thm:ratio.inc}:}
Using (\ref{vmdiffeq}) one finds that
\begin{eqnarray*}\label{diffeq.ratio}
\dffx R_{m+1}(x)=2x \left[ \frac{R_{m+1}(x)}{R_m(x)}-1 \right].
\end{eqnarray*}
After rewriting this in terms of $V_m$ and then using the recursion relation
(\ref{vmprecurs}) with $p=2$ to eliminate
$V_{m+1}$, one finds that $R_{m+1}^{\prime}(x) \geq 0$ if and only if
\begin{eqnarray*}
2(m+1)[R_m(x)]^2-(2m+1-2x^2)R_m(x)-2x^2 \leq 0.
\end{eqnarray*}
The polynomial
$P(z)=2(m+1)z^2-(2m+1-2x^2)z-2x^2 $
has one positive and one negative root, and
$R_{m+1}^{\prime}(x) \geq 0$ if and only if  $R_m(x)$
lies between these two roots.
Since $1 \geq R_m(x) > 0$, it follows that $R_{m+1}(x)$ is increasing if
and only if $R_m(x)$ is less than the larger root, i.e.,
$$R_m(x) \leq
\frac{4x^2}{ \sqrt{4(x^2+m)^2+1+4m+12x^2}+2x^2-2m-1}.$$
where, we have again written the root in the non-standard form
$\frac{C}{-B+\sqrt{B^2 - AC}}$. Then using the upper bound of Theorem
\ref{thm:optm.ineq}, we see that it suffices to show that
\begin{eqnarray*}
R_m(x) & \leq & \frac{4x^2}{\sqrt{(x^2+m)^2+4x^2}+3x^2-m} \\
 & \leq &
\frac{4x^2}{\sqrt{4(x^2+m)^2+1+4m+12x^2}+2x^2-2m-1}
\end{eqnarray*}
or, equivalently that
\begin{eqnarray*}
\lefteqn{ \sqrt{(x^2+m)^2+4x^2}+3x^2-m } ~~~~~~~~~~~~ \\
  & \geq & \sqrt{4(x^2+m)^2+1+4m+12x^2}+2x^2-2m-1
\end{eqnarray*}
which is easily checked.

\section{Proof of Ratio Bounds}\label{sect:rat.bnd.pfs}

The proofs in this section, although elementary, are quite long
and tedious.  The details were checked using Mathematica.

\subsection{Differential inequality}

In order to prove Theorem \ref{thm:optm.ineq}, it suffices
to establish the following

\begin{lemma}\label{ineq:diffGm}  Let $G_k^m$ be given by (\ref{def:Gkm}).
Then
\begin{itemize}
\item[{\rm (i)}] For $m \geq 1 $,
 $\ds{~~~\dffx G_4^m(x^2) \leq 2x(\frac{G_4^m(x^2)}{G_4^{m-1}(x^2)}-1)}$.
\item[{\rm (ii)}] For $m \geq 4 $,
 $\ds{~~~ \dffx G_8^m(x^2) \geq 2x(\frac{G_8^m(x^2)}{G_8^{m-1}(x^2)}-1)},$
\end{itemize}
but the inequality (ii) does not hold  for $m < 4$.
\end{lemma}
\pf  The proof is based on the elementary principle that
if a function on the half-line is zero at the origin
 and increasing, then it is non-negative.
Unfortunately, the actual verification is rather tedious
and requires the repeated use of this principle.
For simplicity, we put $x^2=y$ and assume $y \geq 0$.
Then (i) is equivalent to
\begin{eqnarray}
  E_m(y) \equiv \left( \frac{G_4^m(y)}{G_4^{m-1}(y)}-1 \right) -
   \dffy G_4^m(y) \geq 0.
\end{eqnarray}
Let $B_m=(m^2+y^2+4y+2my)^\frac{1}{2}$.  Then
$$G_4^m(y)=\frac{4y}{B_m+3y-m}. $$
and
$$E_m(y)=\frac{B_m \left[4m+(B_m+3y-m)(B_{m-1}-B_m+1)\right]
-(4m^2+8y+4my)}{B_m~(B_m+3y-m)^2}.$$
Thus $E_m(y) \geq 0 $ if and only if
\begin{eqnarray*}
\lefteqn{B_m\biggr(3m+B_{m-1} [B_m+3y-m]\biggr)+m^3+2my } \\
& \geq & 3m^2+4y+11y^2+m^2y+5my^2+3y^3+B_m\biggl((y+m)^2+y\biggr).
\end{eqnarray*}
We put
\begin{eqnarray*}
 s & = & s(y,m)=3m^2+4y+11y^2+m^2y+5my^2+3y^3, \\
t & = &t(y,m)=2my+m^3,  ~~\hbox{and} \\
h& = &h(y,m)=(y+m)^2+y-3m.
\end{eqnarray*}
Then $E_m(y) \geq 0 $ if and only if
\begin{eqnarray}
B_mB_{m-1}(B_m+3y-m)\geq B_m h + s-t.\label{bm1}
\end{eqnarray}
Notice that both sides of (\ref{bm1}) are positive. For the left
side this follows immediately from $B_m > m$. For the right, note
that $B_m h(0)+s(0)-t(0)=0$ and
\begin{eqnarray*}
\lefteqn{ \dffy  \left[ B_mh(y)+s(y)t(y)\right]} \\
 & = & \frac{1}{B_m} \left[
6y +12my+9m^2y +12y^2+9my^2+3y^3 \right.\\
& ~ & \left. + ~ B_m(22y+10my+9y^2) +
4B_m-2mB_m+3m^3-6m +m^2B_m \right] .
\end{eqnarray*}
Now observe that
\begin{eqnarray*}
\lefteqn{4B_m-2mB_m+3m^3-6m +m^2B_m}~~~~~~~~~ \\
 & = & 3m(m^2-2)+B_m(m^2-2m+4) \geq 3m(m^2-1)  \geq 0.
\end{eqnarray*}
since $m \geq 1$ and $B_m \geq m$.
Hence $B_m h+s-t$ is increasing in $y$ and the right  side of
(\ref{bm1}) is also positive.
Therefore we can square both sides of (\ref{bm1}) to conclude that
it is equivalent to
\begin{eqnarray}
F(y)=B_m f_1(y) -f_2(y) \geq 0, \label{bm2}
\end{eqnarray}
where
\begin{eqnarray*} f_1(y)=(m+y)(y^3+my^2+3y^2-m^2y-3my+2y-m^3+2m^2)
\end{eqnarray*}
and
\begin{eqnarray*}
\lefteqn{f_2(y)= y^5+3my^4+5y^4+2m^2y^3+5my^3 +2y^3} ~~~~~ \\
& ~ & ~~ - ~ 2m^3y^2-3m^2y^2-3m^4y-   m^3y+6m^2y-m^5+2m^4.
\end{eqnarray*}
Note that $F(0)=0$. Therefore, to prove (\ref{bm2}) it is enough
to show that
$\dffy F(y) = B_m f_1^{\prime}(y)
+\frac{f_1(y)(2+y+m)}{B_m}-f_2^{\prime}(y)
\geq 0$, or equivalently
\begin{eqnarray}
D(y) \equiv d_1(y)-B_m d_2(y) \geq 0, \label{bm3}
\end{eqnarray}
where
\begin{eqnarray*}
 d_1(y) & = &  B_m^2 f_1^{\prime}(y) +f_1(y)(2+y+m) \\
& = & 6 m^3 - m^4 - 3 m^5 + 12 m y + 4m^2 y - 13 m^3 y - 7 m^4 y + 20 y^2
\\ & ~ &
+ ~  14 m y^2 + 7 m^2 y^2 + 2 m^3 y^2 +
48 y^3 + 49 m y^3 + 18 m^2 y^3 + 30 y^4
\\ & ~ &   + ~  17 m y^4 + 5 y^5, ~~~\hbox{and} \\
 d_2(y) & = & f_2^{\prime}(y) = 6 m^2 - m^3 - 3 m^4 - 6 m^2 y - ~4 m^3 y
  \\  & ~ &  ~~~~~  + ~ 6 y^2 +
    15 m y^2 + 6 m^2 y^2 + 20 y^3 + 12 m y^3 + 5 y^4.
\end{eqnarray*}
Note that $D(0)=0$. Therefore, to prove (\ref{bm3}) it is enough
to show that
$\dffy D(y) = d_1^{\prime}(y)-B_m d_2^{\prime}(y)
-\frac{d_2(y)(2+y+m)}{B_m}
\geq 0$, or equivalently,
\begin{eqnarray}
G(y) \equiv B_mg_1(y)-g_2(y) \geq 0, \label{bm4}
\end{eqnarray}
where
\begin{eqnarray*}
 g_1(y)& = & d_1^{\prime}(y)  \\
& = & 12 m + 4 m^2 - 13 m^3 - 7 m^4 + 40 y + 28 m y + 14 m^2 y +
    4 m^3 y  \\
& ~ & ~~~ + ~ 144 y^2 + 147 m y^2 + 54 m^2 y^2 + 120 y^3 +
    68 m y^3 + 25 y^4, ~~\hbox{and} \\
g_2(y)& = & B_m^2 d_2^{\prime}(y) +d_2(y)(2+y+m)  \\
& = & 12 m^2 + 4 m^3 - 13 m^4 - 7 m^5 - 18 m^2 y - 13 m^3 y -
    3 m^4 y  \\
 & ~ & + ~  60 y^2 +  180 m y^2 + 183 m^2 y^2 + 58 m^3 y^2 +
    298 y^3 \\  & ~ &
  + ~  353 m y^3 + 122 m^2 y^3 + 170 y^4 + 93 m y^4  +  25 y^5.
\end{eqnarray*}
Note that $G(0)=0$. Therefore, to prove (\ref{bm4}) it is enough
to show that
$\dffy G(y) = B_m g_1^{\prime}(y)
+\frac{g_1(y)(2+y+m)}{B_m}-g_2^{\prime}(y)
\geq 0$, or equivalently
\begin{eqnarray}
H(y)=h_1(y)-B_m h_2(y) \geq 0, \label{bm5}
\end{eqnarray}
where
\begin{eqnarray*}
h_1(y) & = & B_m^2 g_1^{\prime}(y) +g_1(y)(2+y+m)  \\
& = & 24 m + 60 m^2 + 6 m^3 - 13 m^4 - 3 m^5 + 240 y + 300 m y +  460m^2 y
\\  & ~&  ~ + ~
347 m^3 y + 113 m^4 y + 1520 y^2 + 2246 m y^2 + 1663 m^2 y^2 + 482 m^3 y^2
\\ & ~&  ~ + ~
2112 y^3 + 2233 m y^3 +  738 m^2 y^3 + 930 y^4 + 497 m y^4 + 125 y^5,
  ~~\hbox{and} \\
h_2(y)& = & g_2^{\prime}(y)\\
& = & -18 m^2 - 13 m^3 - 3 m^4 + 120 y + 360 m y + 366 m^2 y +
116 m^3 y   \\
& ~&  ~ + ~  894 y^2 +1059 m y^2 + 366 m^2 y^2 + 680 y^3 + 372 m y^3 + 125 y^4.
\end{eqnarray*}
Note that $H(0)=12m(2+5m+2m^2) >0$. Therefore, to prove (\ref{bm5}) it is
enough
to show that
$\dffy H(y) = h_1^{\prime}(y)-B_m h_2^{\prime}(y)
-\frac{h_2(y)(2+y+m)}{B_m}
\geq 0$, or equivalently
\begin{eqnarray}
l_1(y)B_m - l_2(y) \geq 0, \label{bm6}
\end{eqnarray}
where
\begin{eqnarray*}l_1(y) & = & h_1^{\prime}(y)  \\ & = &
240 + 300 m + 460 m^2 + 347 m^3 + 113 m^4 + 3040 y + 4492 m y
\\ &~& + ~ 3326 m^2 y + 964 m^3 y + 6336 y^2 + 6699 my^2 +  2214 m^2 y^2 +
\\ &~& +~  3720 y^3 + 1988 m y^3 + 625 y^4,  ~~~\hbox{and} \\
l_2(y) & = & B_m^2 h_2^{\prime}(y)+h_2(y)(2+y+m) \\
& = & 84 m^2 + 316 m^3 + 347 m^4 + 113 m^5 + 720 y + 2520 m y + 5046 m^2 y
\\ &~& +~
3899 m^3 y + 1077 m^4 y + 9180 y^2 + 15780 m y^2 + 11727 m^2 y^2
\\ &~& +~  3178 m^3 y^2 +
12202 y^3 + 13145 m y^3 + 4202 m^2 y^3 + 4970 y^4
\\ &~& +~ 2613 m y^4 +  625 y^5.
\end{eqnarray*}
Note that $l_1(y) \geq 0$ and $l_2(y) \geq 0$ for all $y \geq 0$.
Therefore (\ref{bm6}) holds, if and only if
$L(y)=B_m^2(l_1(y))^2-(l_2(y))^2 \geq 0$,
which follows immediately from the fact that all the coefficients
are positive in
\begin{eqnarray*}
L(y) & = &
4(14400 m^2 + 36000 m^3 + 75936 m^4 + 97368 m^5 + 78972 m^6 + 37188 m^7 \\
&~& +~ 8136 m^8 + 57600 y + 172800 m y + 717360 m^2 y + 1373400 m^3 y  \\
&~& +~ 1732428
m^4 y +  1360314 m^5 y + 599454 m^6 y + 111444 m^7 y + 1344000 y^2 \\
 &~&+~
3838560 m y^2 +   8437260 m^2 y^2 + 11062920 m^3 y^2 + 8495031 m^4 y^2
\\ &~& +~ 3499083 m^5 y^2  +
595986 m^6 y^2 + 9342880 y^3 + 24217360 m y^3 \\
&~& +~ 32546720 m^2 y^3
+ 24561680 m^3y^3 +  9950080 m^4 y^3 + 1694280 m^5 y^3 \\
 &~& +~ 17918380 y^4
 + 37038224 my^4 + 34271234 m^2 y^4 + 15627870 m^3 y^4 \\
 &~& +~ 2862630 m^4 y^4
+ 15343236 y^5 + 23700982 m y^5 + 13930330 m^2 y^5 \\
 &~& +~ 2982516 m^3 y^5  + 6445963 y^6 + 6618363 m y^6 + 1887294 m^2 y^6\\
&~& +~ 1302640 y^7  + 667056 m y^7 +101250 y^8).
\end{eqnarray*}
To prove (ii) we proceed similarly, but now let
 $B_m(y) =\sqrt{(y+m)^2+8y }$ and
$E_m(y)=\dffy G_8^m(y)-(\frac{G_8^m(y)}{G_8^{m-1}(y)}-1)$,
noting that
$$G_8^m(y)=\frac{8y}{B_m+7y-m} .$$
We now need to show that
$E_m(y) \geq 0$ for all $m \geq 4$.
As the argument is similar to that above, we omit the details
except to  indicate the steps leading to  the condition
$m \geq 4$.
Observe that $E_m(y) \geq 0 $ if and only if
\begin{eqnarray}
\lefteqn{B_m((y+m)^2+y-7m)} \nonumber \\ &  +7 m^2 - m^3 + 24 y -
2 m y + 5 m^2 y + 55 y^2 + 13 m y^2 + 7y^3 \nonumber \\ & \geq B_m
B_{m-1}(7y-m+B_m).\label{bm7}
\end{eqnarray}
Again, both sides of the inequality are positive. Hence we can
square both sides of the inequality and, as above get that
(\ref{bm7}) is equivalent to
\begin{eqnarray}
F(y)=f_1(y)-B_m f_2(y) \geq 0, \label{bm8}
\end{eqnarray}
with the appropriate $f_1$ and $f_2$. Again $F(0)=0$.
Therefore, in order to prove (\ref{bm8}), it is enough to show that
$\dffy F(y) \geq 0$, or equivalently, after rewriting,
\begin{eqnarray}
D(y)=B_m d_1(y)-d_2(y) \geq 0, \label{bm9}
\end{eqnarray}
with the appropriate $d_1$ and $d_2$. And again $D(0)=0$.
We repeat the procedure: To prove (\ref{bm9}), it is enough to show that
$\dffy D(y) \geq 0$, or equivalently,
$$e_1(y)-B_m e_2(y)\geq 0, $$
with the appropriate $e_1$ and $e_2$. $e_1$ and $e_2$
turn out to be  both positive for $y \geq 0$. Therefore (\ref{bm9}) holds
if
$$L(y)=(e_1(y))^2-(B_m e_2(y))^2 \geq 0.$$
 $L(0)=0$ and
$$L^{\prime}(0)=192m^2(m-4)(1+2m)(480+64m+90m^2+33m^3).$$
Thus $L^{\prime}(0) \geq 0 $ if and only if $m \geq 4$. For
all $m \geq 4$, $L^{\prime \prime}(y) \geq 0$ for all $y \geq 0$.
This finishes (ii).   \qed

\subsection{Proof of Theorem \ref{thm:optm.ineq}}

We will prove by induction that $R_m(x) < G_4^m(x^2)$ for $m = 0, 1, 2, 3,
\ldots$
As observed earlier,
this inequality  holds for $m = 0$, since
it is then equivalent to the upper bound in (\ref{opt.ineq}).
Let
$$H_m(x)=G_4^m(x^2)-R_m(x).$$
Then the upper bound in Theorem \ref{thm:optm.ineq} is equivalent to
$H_m(x) \geq 0$.
This can be verified using the strategy of Section \ref{sect:V0} if the
following
conditions hold
\begin{itemize}
  \item [(i)] $H_m(0) >0.$
  \item [(ii)] $\ds{\lim_{x \rightarrow \infty} H_m(x)=0}.$
   \item [(iii)] $H_m^{\prime}(x) \leq F_{\rm pos}(x) H_m(x)$ for some
strictly positive function $F_{\rm pos}(x) > 0$.
\end{itemize}
Conditions (i) and (ii) hold. Indeed,
$$H_m(0)=\frac{2m}{1+2m}-
\frac{\Gamma(m) \Gamma(m+\frac{1}{2})}{\Gamma(m+1)
\Gamma(m-\frac{1}{2})}=\frac{1}{2m(1+2m)},$$
and
$$\ds{\lim_{x \rightarrow \infty} H_m(x)=0},$$
since $\ds{\lim_{x \rightarrow \infty}R_m(x)=1}$, and for all $k \geq 1$,
$\ds{\lim_{x \rightarrow \infty} G_k^m(x)=1}.$

\bigskip
We check now condition (iii).
It follows from Lemma \ref{ineq:diffGm} (i)  and (\ref{diffeq.ratio}) that
\begin{eqnarray*}
H_m^{\prime}(x) & \leq & 2x \left[
  \frac{G_4^m(x^2)}{G_4^{m-1}(x^2)}  - \frac{R_m(x)}{R_{m-1}(x)} \right] \\
   & = & \frac{2x}{G_4^{m-1}(x^2) R_{m-1}(x)} \left[
   G_4^m(x^2) R_{m-1}(x) -G_4^{m-1}(x^2) R_m(x) \right] \\
   & \leq & \frac{2x}{R^{m-1}(x)} \left[ G_4^m(x^2)  - R_m(x) \right] \\
   & = &  \frac{2x}{R^{m-1}(x)} H_m(x)
\end{eqnarray*}
where the inequality follows from the induction hypothesis
$R_{m-1}(x) < G_4^{m-1}(x^2)$.  Thus (iii) holds with
$F_{\rm pos}(x) = 2x/R_{m-1}(x).$

For $m \geq 4$, the lower bound is proved similarly.
One considers $H_m(x)=R_m(x)-G_8^m(x^2)$ instead and uses
Lemma \ref{ineq:diffGm}(ii).
The cases $R_1$, $R_2$ and $R_3$ have to be verified directly
as Lemma \ref{ineq:diffGm} (ii) only covers the cases $R_m$ for $m\geq 4$.

Using (\ref{vmprecurs}), $R_1 \geq G_8^0$ is equivalent to showing that
$$\frac{x}{V_0(x)} \geq \frac{9x+14x^3+(2x^2-1)\sqrt{8+x^2}}
{2(7x+\sqrt{8+x^2})}.$$
As $V_0$ is always positive, this inequality holds trivially for those $x$, for
which the right hand side is negative or zero.
Therefore we only need to prove the inequality on the interval
$[x_0, \infty)$, $x_0 \simeq 0.2511$, where
$$9x+14x^3+(2x^2-1)\sqrt{8+x^2} \geq 0.$$
Hence we need to show
that for all $x \in [x_0, \infty)$,
\begin{eqnarray} \nonumber
V_0(x) \leq 2x
\frac{7x+\sqrt{8+x^2}}{9x+14x^3+(2x^2-1)\sqrt{8+x^2}}\\
=2x\frac{6x^2-1}{1+6x^2+12x^4-2x\sqrt{8+x^2}}.\label{r1}
\end{eqnarray}
Put $h_1(x)=2x\frac{6x^2-1}{1+6x^2+12x^4-2x\sqrt{8+x^2}}.$
By Theorem \ref{ineqv0} of Section \ref{sect:V0},
 inequality (\ref{r1}) is true for all
$x \in[x_0,\infty)$, for which
$$g_4(x) \leq  h_1(x).$$
This last inequality holds only on an interval $[x_0,x_1]$,
$x_1 \simeq 1.399$. For all $x \geq x_1$, we show that
$$h_1(x) < \frac{1}{x}$$
and
$$h_1^{\prime} \leq 2(x h_1 -1).$$
Then (\ref{r1}) follows as in  Section \ref{sect:V0}.

Next, we find that $R_2 \geq G_8^1$ is equivalent to
$V_0(x) \geq h_2(x)$, where
$$h_2(x)=2x \frac{3 + 9  x^2 + 14  x^4 +
      (2  x^2 - 3)(8 x^2 + (1 + x^2)^2))^\frac{1}{2}}
{-3 - 7x^2 + 32 x^4 + 28 x^6 + ( 3 -  4 x^2+ 4 x^4)
(8 x^2 + (1 + x^2)^2))^\frac{1}{2}}$$
and  $R_3 \geq G_8^2$ is equivalent to
$V_0(x) \leq h_3(x)$ where
\begin{eqnarray*}
  h_3(x) =
 \frac{2x}{N(x)} \left[-30 - 23 x^2 + 32 x^4 + 28 x^6 +
  \sqrt{8 x^2 + (2 + x^2)^2}~
 (15 - 8x^2 + 4 x^4) \right],
\end{eqnarray*}
with
\begin{eqnarray*}
\lefteqn{ N(x)=30 + 3 x^2 - 42 x^4 + 92 x^6 + 56 x^8  } \\
& & +  (8 x^2 + (2 + x^2)^2))^\frac{1}{2}(-15 + 18x^2 -12 x^4 +8x^6).
\end{eqnarray*}
Again, we have to check these inequalities for $V_0$  only for those $x$,
for which
the  right hand sides are positive. We then proceed as
for
$R_1$ and show that $g_{\pi} \geq h_2$ up to a certain $x_2$ and
that $h_2 <\frac{1}{x}$, $h_2^{\prime} \geq 2(x h_1 -1)$ on $[x_2, \infty)$.
Similarly, we show that
$g_4 \leq h_3$  up to a certain $x_3$ and
that $h_3 <\frac{1}{x}$, $h_3^{\prime} \leq 2(x h_1 -1)$ on $[x_3, \infty)$.
\qed

Note that these arguments also show that on the interval $[x_1, \infty)$
the function $h_1$ is a better upper bound for $V_0$ than $g_4$;
  on $[x_2, \infty)$ the function $h_2$ is a better lower bound
for $V_0$ than $g_{\pi}$; and on
 $[x_3, \infty)$ the function $h_3$ is a better upper bound
for $V_0$ than $g_4$. In fact,  $h_3 \leq h_1 \leq g_4$ for $x > x_3$.

\subsection{Optimality of bounds}\label{sect:opt}

We still need to consider optimality of the lower bound
in upper bound in (\ref{eq:rat.bnd} in the parameter $k$.
We continue the strategy above
using similar notation so that now
 $B_m=\sqrt{(y+m)^2+ky}$  and
$E_m(y)= \left[\frac{G_k^m(y)}{G_k^{m-1}(y)}-1 \right] -~\dffy
G_k^m(y)$ with $$G_k^m(y)=\frac{ky}{B_m+(k-1)y-m} .$$ Then $E_m(y)
\leq 0$ if and only if,
\begin{eqnarray}\label{bk1}
2B_mB_{m-1}(B_m+(k-1)y-m)\leq  2B_m (y+(y+m)^2-(k-1)m) +P
\end{eqnarray}
where
\begin{eqnarray*}
\lefteqn{P=-2 m^2 + 2 k m^2 - 2 m^3 - 2 k y + k^2 y - 4 m y -
    6 m^2 y + 2 k m^2 y}~~~ \\
& ~ & ~~~~ - ~ 2 y^2 - 2 k y^2 +  2 k^2 y^2 -
    6 m y^2 + 4 k m y^2 - 2 y^3 + 2 k y^3.
\end{eqnarray*}
For $m\geq 1$ and $k \geq 2$  both sides of the inequality are positive.
Therefore we can square both sides and get that $E_m(y) \leq 0$ if and only if
$$F(y)=f_1(y)-B_m f_2(y) \geq 0,$$
where
\begin{eqnarray*}
f_1(y) & = & -16 k m^4 + 8 k^2 m^4 + 32 m^5 - 16 k m^5 -
    20 k^2 m^2 y + 8 k^3 m^2 y + 8 k m^3 y
\\ & ~ & - ~ 4 k^2 m^3 y +
    160 m^4 y - 96 k m^4 y + 8 k^2 m^4 y - 4 k^3 y^2 + k^4 y^2
\\ & ~ & + ~ 120 k m^2 y^2  - 84 k^2 m^2 y^2 +
    12 k^3 m^2 y^2 + 320 m^3 y^2 - 224 k m^3 y^2
\\ & ~ & + ~  32 k^2 m^3 y^2 + 20 k^2 y^3 - 20 k^3 y^3  + 4 k^4 y^3 +
    152 k m y^3 - 124 k^2 m y^3
\\ & ~ & + ~ 24 k^3 m y^3 + 320 m^2 y^3 -
 256 k m^2 y^3 + 48 k^2 m^2 y^3  +
    56 k y^4 - 52 k^2 y^4
\\ & ~ & +~   12 k^3 y^4 + 160 m y^4 - 144 k m y^4
 + 32 k^2 m y^4 + 32 y^5 - 32 k y^5
+  8 k^2y^5,
 \\ f_2(y) & = & 4 \left(
      -4 k m^3 + 2 k^2 m^3 + 8 m^4 - 4 k m^4 -
        3 k^2 m y + k^3 m y + 2 k m^2 y  \right.
\\ & ~ & - ~ k^2 m^2 y + 32 m^3 y - 20 k m^3 y
         + 2 k^2 m^3 y + k^2 y^2 + 16 k m y^2 - 12 k^2 m y^2
\\ & ~ & +~  2 k^3 m y^2 +
        48 m^2 y^2 - 36 k m^2 y^2 + 6 k^2 m^2 y^2 +
        10 k y^3 -  9 k^2 y^3 + 2 k^3 y^3
\\ & ~ & + ~ \left.  32 m y^3 -
       28 k m y^3 + 6 k^2 m y^3 + 8 y^4 - 8 k y^4 +
        2 k^2 y^4\right).
\end{eqnarray*}
Then $F(0)=0$ and in order that $F \geq 0$, we must have $\dffy
F(0) \geq 0$. Computing $\dffy F(y)$, we find that $\dffy F(0) =
0$. We apply the same procedure as in  the proof of Lemma
\ref{ineq:diffGm}, compute the successive derivatives and evaluate
them at $0$. Evaluating the derivative at $0$, in the fourth step
of the procedure gives the value $$24k^3m(1+2m)(km-6m-k).$$
Therefore, in order that (\ref{bk1}) (which is the condition for
the lower bound) holds for all $m\geq 2$,   we have to have at
least that $k\geq \lim_{m \rightarrow \infty} \frac{6m}{m-1} = 6$.
Thus for $m=2$, $k\geq 12$ will do, for $m=3$, $k\geq 9$, for
$m=4$, $k\geq 8$ and so fourth. Therefore, as $G_k^m$ is
increasing in $k$, it seems a natural choice to pick $k=12$ or
bigger for the lower bound. And indeed, one can check that $G_k^m$
satisfies the lower bound condition of Lemma \ref{ineq:diffGm} for
$k \geq 12$ and $m \geq 2$. However, it is not true that for $k >
8$, $G_k^{m-1}$ is a lower bound for $R_m$, for {\em all} $m \geq
1$. It is  a lower bound for all $m \geq m(k)$, from a certain
$m(k)$ on. Thus the induction in the proof of  Theorem
\ref{thm:optm.ineq} cannot start at $m=0$ or $m=1$. For $m <
m(k)$, there exists $x_m$ such that  $R_m- G_k^{m-1} \geq 0$ on
$[0,x_m]$ and $R_m- G_k^{m-1} < 0$ on $(x_m, \infty).$

\subsection{Extensions to  general $p$}\label{sect:genp}

For $p=1$,  all the functions involved are
identically equal to $1$ and hence trivially convex.
For large $x$, $1/V_m^p(x)\approx x^{p-1}$ and
$x^{p-1}$ is concave for $1 < p < 2$.   Hence we can not expect convexity
of $1/V_m^p$  on $(0,\infty)$ for $p$  in $(1,2)$.
It was shown in \cite{Mas} that
$\frac{1}{V_0^p}$ is not convex on ${\bf R ^+}$ for $0 < p < 1$.
Therefore, we study only generalizations to $p > 2$. Our method of proof
yields verification of  the convexity of $\frac{1}{V_m^p}$ for all $m \geq 1$
 up to at least $p=4$. However this method breaks down for
larger $p$.

We generalize our previous notation to
$R_m^p(x)=\frac{V_m^p(x)}{V_{m-1}^p(x)}$, and observe that
\begin{eqnarray}
 \dffx R_m^p(x)
=px^{p-1}\left[\frac{R_m^p(x)}{R_{m-1}^p(x)}-1\right].\label{ineq.inp0}
\end{eqnarray}
For $k \geq 1$, $m \geq 0$, $p\geq 1$ we generalize $G_k^{m}$ to
$$G_k^{m,p}(x^p)=\frac{kpx^p}{p[(k-1)x^p-m] +
\sqrt{ p^2(x^p+m)^2+2kp(p-1)x^p }}.$$

The proofs of Theorems \ref{thm:optm.ineq},  \ref{thm:ratio.inc} and
\ref{thm:convex.m} can be extended provided that the analogue of
Lemma  \ref{ineq:diffGm}  holds. This is not the case for large $p$.
However, although the generalization of the  upper bound in
Theorem \ref{thm:optm.ineq} is a necessary and sufficient
condition for convexity of $1/V_m^p(x)$, Lemma \ref{ineq:diffGm}
is only a sufficient condition for Theorem \ref{thm:optm.ineq}.
Indeed, we were able to establish the lower bound in
Theorem \ref{thm:optm.ineq} for $m=1,2,3$ even though part (ii)
of Lemma \ref{ineq:diffGm} does not hold for $m < 4$.
Hence, the fact that Lemma  \ref{ineq:diffGm} breaks down for
large $p$ does not preclude convexity of $1/V_m^p(x)$.
On the contrary,  numerical evidence
suggests that $\frac{1}{V_m^p}$ is convex for all $p \geq 2$.

\begin{lemma}\label{breakdown}
For all $4 \geq p \geq 2$, $m
\geq 1$,
$$\dffx G_4^{m,p}(x^p) \leq px^{p-1}\left[\frac{G_4^{m,p}(x^p)}
{G_4^{m-1,p}(x^p)}-1\right]$$
\end{lemma}
This is  equivalent to
$$E_m^p=px^{p-1}
\left[\frac{G_4^{m,p}(x^p)}{G_4^{m-1,p}(x^p)}-1\right] -
\dffx G_4^{m,p}(x^p) \geq 0$$
which allows us to make some remarks about the range of
validity.  Although Lemma \ref{breakdown}
can probably be extended to some
higher $p$, it does {\em not} hold for all $p,m$.  On the
contrary, for all $m \geq 1$ there exists
$p(m)$ and an interval $(x_1^{p(m)},x_2^{p(m)})$ such that
$E_m^p < 0$ on that interval for all $p \geq p(m)$.
For example,
numerical results show that for $m =1$, an interval on which
$E_m < 0$ exists when $p \geq 10$; for $m=2$ when $p \geq 14$
and for $m=3$ when $p \geq 18$.

For simplicity, we only sketch the proof of Lemma \ref{breakdown}
and give the final expressions for $p=3$. Similar expressions can
be given for $p=4$. We have checked the details using Mathematica,
but omit the long formulas. As the expressions involved are
monotone in $p$, this suffices for the entire interval $2 \leq p
\leq 4$.

\pf Let
$B_m = \sqrt{p^2(y+m)^2+8p(p-1)y}.$
Then $E_m^p \geq 0$ is equivalent to
\begin{eqnarray*}
\lefteqn{B_mB_{m-1}(B_{m}+3py-pm)}~\\
& \geq & pB_{m}\left[p(y+m)^2+y(5p-8)-3pm \right] +
   p^2 \left\{3m^2p-m^3p \right.\\& ~ & ~+ ~ \left. y(-8+8m+ 8p-6mp+m^2p) +
 y^2(-24+23p+5mp)+3py^3 \right\}.
\end{eqnarray*}
Following the procedure used to prove Lemma \ref{ineq:diffGm} (i),
we eventually find that it would suffice to show that
\begin{eqnarray}
l_1(y)B_m - l_2(y) \geq 0, \label{bmp1}
\end{eqnarray}
 where, for $p=3$,
\begin{eqnarray*} l_1(y) & = & 3( 5120 +
       2880 m + 13800 m^2 + 11034 m^3 + 3051 m^4 + 169600 y  \\
& ~ &  + ~
        199632 m y + 116820 m^2 y + 26028 m^3 y + 298296 y^2 +
        239706 m y^2  \\
& ~ & + ~ 59778 m^2 y^2 + 133920 y^3 + 53676 m y^3 +
        16875 y^4), ~~~~~~\hbox{and} \\
l_2(y)& = & 27 (-2048 m + 192 m^2 + 3448 m^3 + 3678 m^4
    + 1017 m^5 25600 y  \\
& ~ & +~   49920 m y   +  76152 m^2 y + 45330 m^3 y +
        9693 m^4 y + 200000 y^2 \\
& ~ & +~   250152 m y^2 + 139266 m^2 y^2 +
        28602 m^3 y^2 + 195880 y^3 \\
& ~ & +~ 157254 m y^3 +
        37818 m^2 y^3 +  59640 y^4 + 23517 m y^4 + 5625 y^5).
\end{eqnarray*}
Both $l_1(y)$ and $l_2(y)$ in (\ref{bmp1}) are positive, and hence
(\ref{bmp1}) is equivalent to
\begin{eqnarray*}
L=l_1(y)^2B_m^2 - l_2(y)^2 \geq 0. \label{bmp2}
\end{eqnarray*}
which can be verified for  $m \geq 1$. \qed

\bigskip

\noindent {\bf Acknowledgment:}  It is a pleasure to thank Professor
S. Kwapien  for the argument leading to the upper bound in Property (a),
and Professors R. Askey and M. Ismail for helpful discussions
about the properties of hypergeometric functions.
We would also like to thank  Professor V. Mascioni for providing a
copy of \cite{Mas} before publication, and Professor Fink for drawing
our attention to reference \cite{M}.

\end{document}